\documentclass[11pt]{amsart}
\usepackage{amsmath}
\usepackage{amssymb}
\usepackage{amsfonts}
\usepackage{amsthm}
\usepackage{enumerate}
\usepackage[mathscr]{eucal}
\usepackage{verbatim}
\usepackage{amsthm}
\usepackage{amscd}
\usepackage[mathscr]{eucal}
\usepackage{appendix}
\usepackage{tikz}
\usepackage{soul}

\newtheorem{innercustomgeneric}{\customgenericname}
\providecommand{\customgenericname}{}

\newcommand{\newcustomtheorem}[2]{\newenvironment{#1}[1]
  {\renewcommand\customgenericname{#2}
   \renewcommand\theinnercustomgeneric{##1}\innercustomgeneric}{\endinnercustomgeneric}}

\newcustomtheorem{customthm}{Theorem}

\newcommand{\newcustomlemma}[2]{\newenvironment{#1}[1]
  {\renewcommand\customgenericname{#2}
   \renewcommand\theinnercustomgeneric{##1} \innercustomgeneric}{\endinnercustomgeneric}}

\newcustomlemma{customlemma}{Lemma}

\newcustomlemma{customproposition}{Proposition}

\newcommand*{\medcap}{\mathbin{\scalebox{1.5}{\ensuremath{\cap}}}}%

\theoremstyle{plain}
\newtheorem{theorem}{Theorem}

{}
\newtheorem{lemma}{Lemma}\numberwithin{lemma}{section}
\newtheorem{proposition}{Proposition}\numberwithin{proposition}{section}
\numberwithin{corollary}{section}
\theoremstyle{definition}
\numberwithin{definition}{section}
\theoremstyle{remark}
\newtheorem{remark}{Remark}

\numberwithin{equation}{section}
\newcommand{\lmd}{\lambda}
\newcommand{\la}{\lambda}
\newcommand{\ga}{\gamma}
\newcommand{\si}{\sigma}

\newcommand{\R}{\mathbb{R}}
\newcommand{\Z}{\mathbb{Z}}

\newcommand{\cL}{\mathcal{L}}

\newcommand{\q}{\quad}
\newcommand{\qq}{\qquad}

\newcommand{\wh}{\widehat}
\newcommand{\supp}{\mathrm{supp}}

\newcommand{\bbn}{\mathbb{N}}
\newcommand{\bbz}{\mathbb{Z}}

\newcommand{\bbr}{\mathbb{R}}
\newcommand{\bbrn}{\mathbb R^n}

\newcommand{\xxxi}{\vec{\xi}}
\newcommand{\ccdot}{\vec{\cdot}\;}

\def\xxxi{\vec{\boldsymbol{\xi}}}
\def\000{\vec{\boldsymbol{0}}}

\def\kkk{\vec{\boldsymbol{k}}}

\def\yyy{\vec{\boldsymbol{y}}}
\def\zzz{\vec{\boldsymbol{z}}}

\def\ii{{\mathrm{i}}}

\newcommand{\LL}{\mathcal{L}}
\newcommand{\II}{\mathcal{I}}

\newcommand{\Ga}{\Gamma}
\newcommand{\La}{\Lambda}

\title{On the boundedness of multilinear Fourier multipliers on Hardy spaces}

\subjclass[2020]{Primary 42B15, 42B25, 42B30, 47H60}
\keywords{Fourier muitipliers, Multilinear operators, Hardy spaces}

\author{Jin Bong Lee}
\address{J.B. Lee, Research Institute of Mathematics, Seoul National University,Seoul 08826, Republic of Korea}
\email{jinblee@snu.ac.kr}

\author{Bae Jun Park}
\address{B. Park, Department of Mathematics, Sungkyunkwan University, Suwon 16419, Republic of Korea}
\email{bpark43@skku.edu}

\thanks{ J.B. Lee is supported by NRF grant 2021R1C1C2008252. B. Park is supported in part by NRF grant 2019R1F1A1044075 and was supported in part by a KIAS Individual Grant MG070001 at the Korea Institute for Advanced Study } 

\begin{document}

\begin{abstract}
In this paper, we study multilinear Fourier multiplier operators on Hardy spaces. 
In particular, we prove that the multilinear Fourier multiplier operator of H\"ormander type is bounded from $H^{p_1} \times \cdots \times H^{p_m}$ to $H^p$ for $0<p_1,\dots,p_m\le 1$ with $1/p_1 + \cdots 1/p_m = 1/p$, under suitable cancellation conditions. 
As a result, we extend the trilinear estimates in \cite{Lee_Park2021} to general multilinear ones and improve the boundedness result in \cite{LHHLPPY} in limiting situations.
\end{abstract}

\maketitle

\section{Introduction}

Let $m$ be a positive integer greater than or equal to $2$ and $\sigma$ be a bounded function on $(\bbrn)^m$. The $m$-linear Fourier multiplier operator $T_{\sigma}$ associated with $\sigma$ is defined by
\begin{equation}\label{multidef}
T_{\sigma}\big( f_1,\dots,f_m\big)(x)=\int_{(\bbrn)^m}   \sigma(\xi_1,\dots,\xi_m)\wh{f_1}(\xi_1)\cdots \wh{f_m}(\xi_m)e^{2\pi i\langle x,\xi_1+\dots+\xi_m\rangle}    d\xi_1 \cdots d\xi_m
\end{equation}
for Schwartz functions $f_1,\dots,f_m$ on $\bbrn$, where $\wh{f}(\xi):=\int_{\bbrn}f(x)e^{-2\pi i\langle x,\xi\rangle}dx$ is the Fourier transform of $f$.
The  $L^{p_1}\times \cdots \times L^{p_m}\to L^p$ boundedness of $T_{\sigma}$, $1/p=1/p_1+\dots+1/p_m$, is one of principal questions in multilinear harmonic analysis and there have been several satisfactory results even though multilinear operators can not be directly attacked by tools for linear operators.
Coifman and Meyer \cite{Co_Me1978} proved that if $N$ is a sufficiently large number and $\sigma$ satisfies
\begin{equation}\label{multilinearmihlin}
\big| \partial_{\xi_1}^{\alpha_m}\cdots \partial_{\xi_m}^{\alpha_m}\sigma(\xi_1,\dots,\xi_m)\big|\lesssim_{\alpha_1,\dots,\alpha_m} \big| (\xi_1,\dots,\xi_m)\big|^{-(|\alpha_1|+\dots+|\alpha_m|)}, \q (\xi_1,\dots,\xi_m)\not= \000
\end{equation} for all multi-indices $\alpha_1,\dots,\alpha_m$ satisfying $|\alpha_1|+\cdots+ |\alpha_m|\le N$, then $T_{\sigma}$ is bounded from $L^{p_1}(\bbrn)\times \cdots\times L^{p_m}(\bbrn)$ into $L^p(\bbrn)$  for $1<p_1,\dots,p_m<\infty$ and $1\le p<\infty$ with $1/p=1/p_1+\dots+1/p_m$.
This result was refined by Tomita \cite{Tom2010} who obtained a multilinear version of H\"ormander's multiplier theorem.
Let $\Psi$ be a Schwartz function, defined on $(\bbrn)^m$, having the properties that
$0\le \wh{\Psi}\le 1$,  $\supp(\wh{\Psi})\subset \big\{\xxxi\in (\bbrn)^m : 2^{-1}\le |\xxxi|\le 2 \big\}$, and $\sum_{j\in\bbn}\wh{\Psi_j}(\xxxi)=1$ for $\xxxi\not= \vec{0}$
where $\Psi_j:=2^{jmn}\Psi(2^j\ccdot )$. That is, $\{\Psi_j\}_{j\in\bbz}$ is a family of Littlewood-Paley functions on $(\bbrn)^m$.
For $s\ge 0$, we define the Sobolev space $L_s^2((\bbrn)^m)$ in terms of the finiteness of the norm 
\begin{equation}\label{defsobolev}
\Vert F\Vert_{L_s^2((\bbrn)^m)}:=\Big( \int_{(\bbrn)^m}\big(1+4\pi^2|\xxxi|^2 \big)^{s}\big| \wh{F}(\xxxi)\big|^2 \; d\xxxi \Big)^{1/2}
\end{equation}
and then we will use the notation
$$\LL_s^2[\sigma]:=\sup_{j\in\bbz}\big\Vert \sigma(2^j\ccdot)\wh{\Psi}\big\Vert_{L^2_s((\bbrn)^m)}$$
throughout this work.

\begin{customthm}{A}\cite{Tom2010}\label{tomitathm}
Let $1< p, p_1, \cdots, p_m<\infty$ with $1/p = {1}{p_1}+\dots+1/p_m$. Suppose that  
$$s> \frac{mn}{2}.$$
Then we have
\begin{equation*}
\big\Vert  T_\sigma(f_1,\dots,f_m) \big\Vert_{L^p(\R^n)} \lesssim \LL_s^2[\sigma] \prod_{\mu=1}^{m}\Vert f_{\mu}\Vert_{L^{p_{\mu}}(\bbrn)}
\end{equation*}
for Schwartz functions $f_1,\dots,f_m$ on $\bbrn$. 
\end{customthm}
The above result was generalized to the full range of indices $0<p_1,\dots,p_m\le \infty$, in some recent papers, by using a product-type Sobolev space condition, instead of  $\LL_s^2[\sigma]<\infty$. For this one, the space $L^{p_{\mu}}(\bbrn)$, $\mu=1,\dots,m$, would be replaced by the Hardy space $H^{p_{\mu}}(\bbrn)$.\\

To recall the definition of Hardy spaces, 
let $\phi$ be a smooth function defined on $\bbrn$ with $\supp(\phi)\subset \{x\in\bbrn: |x|\le 1\}$. 
Then  the Hardy space
$H^p(\bbrn)$ for $0<p\le \infty$ contains all tempered distributions $f$ on $\bbrn$ such that
\begin{equation}\label{hardydef2}
\Vert f\Vert_{H^p(\bbrn)}:= \Big\Vert \sup_{l\in\Z} \big| \phi_l\ast f\big| \Big\Vert_{L^p(\bbrn)}<\infty
\end{equation} where $\phi_l:=2^{ln}\phi(2^l\cdot)$ for $l\in\bbz$.
The space provides an extension to $0<p\le 1$ in the scale of classical $L^p$ spaces for $1<p\le \infty$, which is more natural and useful in many respects than the corresponding $L^p$ extension.  Indeed, $L^p(\bbrn)=H^p(\bbrn)$ for $1<p\le\infty$ and several essential operators, such as singular integrals of Calder\'on-Zygmund type, that are well-behaved on $L^p(\bbrn)$ only for $1<p\le\infty$ are also well-behaved on $H^p(\bbrn)$ for $0<p\le 1$.
Moreover, especially when $0<p\le 1$,  every $f\in H^p(\bbrn)$ can be written as
\begin{equation}\label{hardydecomp}
f  = \sum_k \lambda_k a_k \qq \text{ in the sense of tempered distributions}
\end{equation}
where $\lambda_k$ are coefficients satisfying $\big(\sum_{k\in\bbn}|\la_k|^p\big)^{1/p}\lesssim \Vert f\Vert_{H^p(\bbrn)}$ and 
$a_k$ are $H^p$-atoms, which means that there exist cubes $Q_k$ such that $\supp(a_k)\subset Q_k$, $\Vert a_k\Vert_{L^{\infty}(\bbrn)}\le |Q_k|^{-1/p}$, and 
\begin{equation}\label{hardyatomvanishing}
\int_{\bbrn}x^{\alpha}a_k(x)dx=0 \qq \text{ for all multi-indices }~ |\alpha|\le L
\end{equation}
 for $L>n/p-n$.
Here, $L$ may be actually an arbitrarily large number greater than $n/p-n$. 
Now let $\mathscr{S}(\bbrn)$ denote the Schwartz space on $\bbrn$ and for a given positive integer $L$ let $\mathscr{S}^{L}(\bbrn)$ be its subspace consisting of $f$ satisfying
$$\int_{\bbrn}x^{\alpha}f(x) dx=0 \q \text{ for all multi-indices}~\alpha~\text{with }~|\alpha|\le L.$$
Then as mentioned in \cite[Proposition 2.1]{Gr_Na_Ng_Sa2019}, for all $0<p<\infty$
\begin{equation}\label{densesubset}
\mathscr{S}^L(\bbrn) ~ \text{ is dense in }~ H^p(\bbrn) \q \text{ if } ~L>\frac{n}{p}-n.
\end{equation}
We remark that $\mathscr{S}(\bbrn)$ is also dense in $H^p(\bbrn)=L^p(\bbrn)$ for $1<p<\infty$, but not for $0<p\le 1$. See \cite[Chapter III, \S 5.2]{St1993} for more details.
Moreover, as discussed in \cite[Chapter III, \S 5.4]{St1993}, if $f\in L^1(\bbrn)\cap H^p(\bbrn)$ for $0<p\le 1$, then we must necessarily have
\begin{equation}\label{hardyvanishing}
\int_{\bbrn}x^{\alpha}f(x) dx=0 \q \text{ for all multi-indices }~|\alpha|\le \frac{n}{p}-n.
\end{equation}
We refer to \cite{Bu_Gu_Si1971, Ca1977, Fe_St1972, St1993, Uch1985} for more details.\\

To study multilinear multiplier operators on Hardy spaces, we now assume that the operator $T_{\sigma}$, defined in \eqref{multidef}, initially acts on functions $f_{\mu}\in \mathscr{S}^L(\bbrn)$, $\mu=1,\dots,m$, in view of \eqref{densesubset}. Here, we set $L$ to be a sufficiently large number greater than $n/p-n$. That is, $L>n/p-n\ge n/p_{\mu}-n$ for all $\mu=1,\dots,m$.
Then we have the following result.
\begin{customthm}{B}\cite{Gr_Mi_Ng_Tom2017, Gr_Mi_Tom2013, Gr_Ng2016, Mi_Tom2013}
Let $0<p<\infty$ and $0<p_1, \dots, p_m\leq \infty$ with ${1}/{p} =1/p_1+\dots+1/p_m$. Suppose that $L$ is sufficiently large and
\begin{equation}\label{prod sob hor condi}
s_1, \dots, s_m > \frac{n}{2}, \qq \sum_{\mu \in J} \Big( s_{\mu} - \frac{n}{p_{\mu}} \Big) > - \frac{n}{2}
\end{equation}
for any nonempty subset $J$ of $ \{1, \cdots, m\}$. 
Then we have
\begin{align}\label{prod sob hor condi esti}
\big\Vert T_{\sigma}(f_1,\dots,f_m)\big\Vert_{L^p(\bbrn)}\lesssim \sup_{j\in \bbz}\big\Vert \sigma(2^j\ccdot)\wh{\Psi} \big\Vert_{L^2_{(s_1,\dots,s_m)}((\bbrn)^m)}\prod_{\mu=1}^{m}\Vert f_{\mu}\Vert_{H^{p_{\mu}}(\bbrn)}
\end{align}
for $f_1,\dots, f_m\in\mathscr{S}^L(\bbrn)$.
\end{customthm}
Here, the space $L_{(s_1,\dots,s_m)}^{2}((\bbrn)^m)$ indicates the product type Sobolev space on $(\bbrn)^m$, in which the norm is defined by replacing the term $(1+4\pi^2 |\xxxi|)^{s}$ in \eqref{defsobolev} by $\prod_{\mu=1}^{m}\big(1+4\pi^2|\xi_{\mu}|^2\big)^{s_{\mu}}$ for $\xxxi:=(\xi_1,\dots, \xi_m)\in (\bbrn)^m$.
We remark that the condition \eqref{prod sob hor condi} is sharp in the sense that if the condition does not hold, then there exists $\sigma$ such that the estimate \eqref{prod sob hor condi esti} fails for the corresponding operators $T_\sigma$.  Refer to \cite{Park_potential} for more details.\\

 In the paper of Lee, Heo, Hong, Park, Yang, and the authors \cite{LHHLPPY}, there has been some recent progress concerning multilinear multiplier theory associated with the original H\"ormander condition, which provides a generalization of Theorem \ref{tomitathm}.
\begin{customthm}{C}\cite{LHHLPPY}\label{ann}
Let $0<p<\infty$ and $0<p_1, \dots, p_m\leq \infty$ with ${1}/{p} = 1/p_1+\dots+1/p_m$. 
Suppose that $L$ is sufficiently large and
\begin{equation}\label{l2 min smooth condi}
s > \frac{mn}{2}, \qq  \frac{1}{p}-\frac{1}{2}<\frac{s}{n}+\sum_{\mu\in J}\Big(\frac{1}{p_{\mu}}-\frac{1}{2}\Big)
\end{equation}
for any subset $J$ of $\{1,\dots,m\}$,  which may be empty, where we adopt the convention that the sum over the empty set is zero.
Then we have
\begin{equation}\label{lp multi esti}
\big\| T_\sigma(f_1, \cdots, f_m)\big\|_{L^p(\R^n)} \lesssim \cL_s^2[\sigma] \, \prod_{\mu =1}^m \| f_{\mu}\|_{H^{p_{\mu}}(\R^n)}
\end{equation}
for $f_1,\cdots,f_m\in \mathscr{S}^L(\R^n)$.
\end{customthm}
The optimality of the condition \eqref{l2 min smooth condi} is already studied by Grafakos, He, and Honz\'ik \cite{Gr_He_Ho2018} who proved that if \eqref{lp multi esti} holds, then we must necessarily have $s\ge mn/2$ and $1/p-1/2\le s/n+\sum_{\mu\in J}(1/p_{\mu}-1/2)$ for all subsets $J$ of $\{1,\dots,m\}$.\\

The estimate \eqref{lp multi esti} has been recently improved in the trilinear case by replacing the target space $L^p(\bbrn)$ by the Hardy space $H^p(\bbrn)$ for $0<p\le 1$, under the additional cancellation condition
\begin{equation}\label{trivanishing}
\int_{\bbrn} x^{\alpha}T_{\sigma}\big(f_1,f_2,f_3\big)(x) dx=0 \qq\text{ for multi-indices }~ |\alpha|\le \frac{n}{p}-n
\end{equation}
for all $f_1,f_2,f_3\in \mathscr{S}^L(\bbrn)$.
\begin{customthm}{D}\cite{Lee_Park2021}\label{trilinear thm}
Let $0<p_1,p_2,p_3<\infty$ and $0<p\le 1$ with $1/p=1/p_1+1/p_2+1/p_3$.
Suppose that $L$ is sufficiently large and
\begin{equation*}
s>\frac{3n}{2} \q \text{ and }\q \frac{1}{p}-\frac{1}{2}<\frac{s}{n}+\sum_{\mu \in J}\Big(\frac{1}{p_{\mu}}-\frac{1}{2} \Big)
\end{equation*}
where $J$ is an arbitrary subset of $\{1,2,3\}$,  which may be empty, where we adopt the convention that the sum over the empty set is zero.
Let $\sigma$ be a function on $(\R^n)^3$ satisfying \eqref{trivanishing}.
Then we have
\begin{equation*}
\big\Vert T_{\sigma}(f_1,f_2,f_3)\big\Vert_{H^{p}(\R^n)}\lesssim \mathcal{L}_s^2[\sigma]\Vert f_1\Vert_{H^{p_1}(\R^n)}\Vert f_2\Vert_{H^{p_2}(\R^n)}\Vert f_3\Vert_{H^{p_3}(\R^n)}
\end{equation*} 
 for $f_1,f_2,f_3\in \mathscr{S}^L(\R^n)$.
\end{customthm}
\begin{remark}
The analogous result is still true in the bilinear setting, but the arguments used there could not be applied to general $m$-linear operators for $m\ge 4$.
\end{remark}
\begin{remark}
The initial function space $\mathscr{S}^L(\bbrn)$ in Theorem \ref{trilinear thm} is actually $\mathscr{S}_0(\bbrn):=\bigcap_{L\in\bbn}\mathscr{S}^L(\bbrn)$ in the original version \cite[Theorem 1]{Lee_Park2021}, but it is easy to see that everything still works with $\mathscr{S}^L(\bbrn)$ as long as $L>n/p-n$. 
In the present work, we employ $\mathscr{S}^L(\bbrn)$ to avoid unnecessary situations.
\end{remark}
\begin{remark}
The condition \eqref{l2 min smooth condi} is equivalent to
$$
s > \frac{mn}{2}, \qq  \frac{1}{p}-\frac{1}{2}<\frac{s}{n}+\sum_{\mu=1}^m \min\Big\{ 0, \Big(\frac{1}{p_{\mu}}-\frac{1}{2}\Big) \Big\}.
$$
Note that $J$ would be empty when $p_1,\dots p_m \leq 2$. \\
\end{remark}

The purpose of this paper is to establish multilinear extensions of Theorem \ref{trilinear thm} in some limiting cases. 
Precisely, we obtain $H^{p_1}\times\cdots\times H^{p_m}\to H^p$ estimates for general $m$-linear operators when $0<p_{1},\dots,p_m \le 1$. 
As mentioned in \cite{Lee_Park2021}, the condition \eqref{trivanishing} is very natural and necessary in the boundedness into the Hardy space $H^p$ for $0<p\le 1$, in view of \eqref{hardyvanishing}.
For a multilinear extension of the above boundedness result, we now consider a bounded function $\sigma$ on $(\bbrn)^m$ having the property
\begin{equation}\label{multivanishing}
\int_{\bbrn} x^{\alpha}T_{\sigma}\big(f_1,\dots,f_m\big)(x) dx=0 \qq\text{ for multi-indices }~ |\alpha|\le \frac{n}{p}-n
\end{equation}
 where $f_1,\dots,f_m\in \mathscr{S}^L(\R^n)$ for $L>n/p-n$.
Then our main result is
\begin{theorem}\label{main}
Let $0< p_1, \cdots, p_m \leq1$ and $0<p\le 1/m$ with $1/p = 1/p_1 + \cdots 1/p_m$. 
Suppose that $L$ is sufficiently large and
\begin{equation}\label{npn2condition}
s>\frac{n}{p} - \frac{n}{2}.
\end{equation}
We also assume that the operator $T_{\sigma}$ satisfies the vanishing moment condition \eqref{multivanishing}. 
Then we have
\begin{equation*}
\big\Vert T_{\sigma}(f_1,\cdots,f_m)\big\Vert_{H^{p}(\R^n)}\lesssim \mathcal{L}_s^2[\sigma]\prod_{\mu=1}^{m}\Vert f_{\mu}\Vert_{H^{p_{\mu}}(\R^n)}
\end{equation*} 
 for $f_1, \cdots ,f_m\in \mathscr{S}^L(\R^n)$.
\end{theorem}
We note that the condition \eqref{npn2condition} clearly implies \eqref{l2 min smooth condi} when $0<p_1,\dots,p_m\le 1$. \\

The study of boundedness into the Hardy space $H^p(\bbrn)$, $0<p\le 1$, for multilinear operators with vanishing integral, such as \eqref{multivanishing}, was initiated by Coifman, Lions, Meyer and Semmes \cite{Co_Li_Me_Se1993} who proved that certain bilinear operators are bounded from $L^{q}\times L^{q'}$ into $H^1$ for $1<q<\infty$ with $1/q+1/q'=1$. Later, the result was generalized by Dobyinski \cite{Do1992} for Coifman-Meyer multiplier operators and by Coifman and Grafakos \cite{Co_Gr1992} for finite sums of products of Calder\'on-Zygmund operators. These works, for $H^1$-boundedness, are related to $compensated$ $compactness$ which arises in theory of partial differential equations to assure convergence of certain nonlinear quantities. For more details, we refer to \cite{Ba_Mu1984, Co_Li_Me_Se1993,  Mul1990, Tar1979} and references therein.
 In the works of Grafakos, Nakamura, Nguyen and Sawano \cite{Gr_Na_Ng_Sa2019, Gr_Na_Ng_Sa_appear}, the $H^{p_1}\times \cdots\times H^{p_m}\to H^p$ boundedness  was established for multilinear multiplier operators $T_{\sigma}$ with some cancellation condition, provided that
$\sigma$ satisfies the multilinear Mihlin's condition \eqref{multilinearmihlin} for multi-indices $\alpha_1,\dots,\alpha_m$ satisfying $|\alpha_1|+\dots+|\alpha_m|\le N$ for sufficiently large $N$. However, the optimal regularity conditions were not pursued at all as it requires sufficiently large $N$.
From this point of view, Theorem \ref{trilinear thm} improves the results of Grafakos, Nakamura, Nguyen, and Sawano for trilinear operators, and Theorem \ref{main} does for general multilinear operators although the input function spaces are limited.\\

The proof of Theorem \ref{main} heavily relies on the atomic decomposition \eqref{hardydecomp},
 but it is not guaranteed that we can interchange infinite sums of atoms and the action of the operator $T_{\sigma}$ as in 
\begin{equation*}
T_{\sigma}\big(f_1,\dots,f_m\big)=\sum_{k_1,\dots,k_m\in\bbn}\la_{1,k_1}\cdots\la_{m,k_m}T_{\sigma}\big(a_{1,k_1},\dots,a_{m,k_m}\big)
\end{equation*}
for functions $f_{\mu}\in H^{p_{\mu}}(\bbrn)$ with atomic representation $f_{\mu}=\sum_{k_{\mu}\in\bbn}{\la_{\mu,k_{\mu}}a_{\mu,k_{\mu}}}$, $1\le \mu\le m$.
 To settle this issue, we will additionally assume that the multiplier $\sigma$ has compact support in $(\bbrn)^m$, which allows us to write, with the aid of Theorem \ref{ann},
$$\big| \phi_l\ast T_{\sigma}\big(f_1,\dots,f_m\big)(x)\big|\le \sum_{k_1,\dots,k_m\in\bbz} \Big(\prod_{\mu=1}^{m}|\la_{\mu,k_{\mu}}|\Big)\big|\phi_l\ast T_{\sigma}\big(a_{1,k_1},\dots,a_{m,k_m}\big)(x) \big|.$$
See Section \ref{reductionsection} below for the details.
Then we apply the method of Grafakos and Kalton \cite{Gr_Ka2001} and pointwise estimates of the form
\begin{equation*}
\big| \phi_l\ast T_{\sigma}\big(a_{1,k_1},\dots,a_{m,k_m}\big)(x)\big|\lesssim \LL_s^2[\sigma]b_1(x)\cdots b_m(x) \q\text{ uniformly in }~ l\in\bbz
\end{equation*} where
$\Vert b_{\mu}\Vert_{L^{p_{\mu}}(\bbrn)}\lesssim 1 $ for $\mu=1,\dots,m$.  
Such arguments are stated in Lemma \ref{keylemma} below, which is actually the key part of this paper.
Since the above estimate separates the left-hand side to $m$ distinct functions of $x$, we may now apply H\"older's inequality with exponents $1/p=1/p_1+\dots+1/p_m$.
For general multipliers $\sigma$ ( without the compact support condition ), we apply a suitable approximation argument, which will be presented in \eqref{generalapprox}.

\subsection*{Notation}
Let $\bbn$ and $\bbz$ denote the sets of the natural numbers and integers, respectively, and we define $\bbn_0:=\bbn\cup \{0\}$. We use the symbol  $A\lesssim B$ to denote $A \leq CB$ for some independent constant $C$, and $A\sim B$ if $A\lesssim B$ and $B\lesssim A$ hold simultaneously.
For a given cube $Q$ in $\bbrn$, let $c_Q$, $\ell(Q)$ denote its center and side length, and $\chi_{Q}$ be the characteristic function of $Q$. 
We denote by $Q^*$, $Q^{**}$, $Q^{***}$ the concentric dilates of a cube $Q$ with $\ell(Q^{*}) = 10\sqrt{n} \ell(Q)$, $\ell(Q^{**})=\big(10\sqrt{n}\big)^2\ell(Q)$, and $\ell(Q^{***})=\big(10\sqrt{n}\big)^3\ell(Q)$.
We use the notation $\langle y\rangle:=\big(1+4\pi^2|y|^2\big)^{1/2}$ for $y\in \bbr^M$, $M\in\bbn$.

\hfill

\section{Preliminaries}\label{prelim}

For a locally integrable function $f$ on $\bbrn$, let
$$\mathcal{M}f(x):=\sup_{Q:x\in Q}\frac{1}{|Q|}\int_Q |f(y)|dy$$
be the Hardy-Littlewood maximal function of $f$ where the supremum ranges over all cubes in $\bbrn$ containing $x$.
Then it turns out that the maximal operator $\mathcal{M}$ is bounded in $L^p(\bbrn)$ for all $1<p\le \infty$. 

We now review some lemmas that will be repeatedly needed in the proof of the main result.
The first lemma is very useful in estimating the $L^p$ norm of the sum of functions having a compact support for $0<p\le 1$. 
\begin{lemma} \cite[Lemma 2.1]{Gr_Ka2001} \label{reverse FS ineq}
Let $0<p\le 1$ and $J$ be a finite or countable index set. 
Suppose that $\{Q_{\mu}\}_{\mu\in J}$ is a sequence of cubes in $\bbrn$, and $\{F_{\mu}\}_{\mu\in J}$ is a family of nonnegative integrable functions with $\supp(F_\mu)\subset Q_{\mu}$.
Then we have
\begin{equation*}
\Big\|\sum_{\mu\in J} \chi_{Q_{\mu}} F_{\mu} \; \Big\|_{L^p(\bbrn)} \lesssim \bigg\| \sum_{\mu\in J} \Big( \frac{1}{|Q_{\mu}|}\int_{Q_{\mu}} F_{\mu}(y)dy \Big) \chi_{Q_{\mu}} \bigg\|_{L^p(\bbrn)}
\end{equation*}
where the constant in the inequality depends only on $p$.
\end{lemma}

We now introduce a variant of Bernstein's inequality, which is contained in \cite{Gr_Mi_Ng_Tom2017, Mi_Tom2013}, modulo a few minor modifications.
\begin{lemma}\cite{Gr_Mi_Ng_Tom2017, Mi_Tom2013} \label{imbedding}
Let $1\le p\le q \le \infty$ and $s\ge 0$. Let $K$ be a function, defined on $(\bbrn)^m$, whose Fourier transform is supported in a ball of a constant radius.
Then for $\mu=1,\dots,m$ and any multi-indices $\alpha\in (\bbn_0)^n$, there exists a constant $C_{p,q,{\alpha}}$ such that
\begin{align*}
\big\|  \big\langle (x_1,\dots,x_m) \big\rangle^{s} \partial_{{\mu}}^{{\alpha}}K(x_1,\dots,x_m) \big\|_{L^q(x_{\mu})} \le C_{p,q,{\alpha}} \big\|  \big\langle (x_1,\dots,x_m) \big\rangle^{s} K(x_1, \cdots, x_m) \big\|_{L^p(x_\mu)}
\end{align*}
 for all  $x_j \in \bbrn, j\not= \mu$.
\end{lemma}

The last one will play a significant role in the proof of Theorem \ref{main}, in conjunction with the vanishing moment condition of $T_{\sigma}$ in \eqref{multivanishing}.
\begin{lemma}\cite[Lemma 4.1]{Lee_Park2021} \label{eps control}
Let $N\in\mathbb{N}$ and $a\in\R^n$.
Suppose that a function $F$, defined on $\bbrn$, satisfies
\begin{equation*}
\int_{\R^n}{x^{\alpha}F(x)}dx=0 \quad \text{for all multi-indices}~ \alpha ~ \text{ with }~|\alpha|\leq N.
\end{equation*}
Then for any $0\leq \epsilon\leq 1$, there exists a constant $C_{\epsilon}>0$ such that
\begin{equation*}
\big\Vert \phi_l \ast F\big\Vert_{L^{\infty}(\R^n)}\le C_{\epsilon} 2^{l(N+n+\epsilon)}\int_{\R^n}{|y-a|^{N+\epsilon}|F(y)|}dy.
\end{equation*}
\end{lemma}

\section{Proof of Theorem \ref{main} : Reduction}\label{reductionsection}

We first show that the multiplier $\sigma$ can be assumed to be compactly supported in Theorem \ref{main}, using an approximation argument.
Suppose that
\begin{equation}\label{conditionsnpn2}
\LL_s^2[\sigma]<\infty \q \text{ for }~ s>n/p-n/2.
\end{equation}
Let $\varphi$ be a Schwartz function on $\bbrn$ such that $\supp(\wh{\varphi})\subset \{\xi\in\bbrn: |\xi|\le 2\}$ and
$\wh{\varphi}(\xi)=1$ for $|\xi|\le 1$. 
For $R>0$, let $\varphi_R:=R^n\varphi(R\, \cdot\,)$ and
$$\wh{\Phi_R}(\xi_1,\dots,\xi_m):=\wh{\varphi_R}(\xi_1)\cdots \wh{\varphi_R}(\xi_m) \q \text{ for }~(\xi_1,\dots,\xi_m)\in (\bbrn)^m.$$
Now we define 
$$\sigma_R(\xxxi):=\sigma(\xxxi)\wh{\Phi_R}(\xxxi)$$
so that
\begin{equation}\label{sigmarrepresentation}
T_{\sigma_R}(f_1,\dots,f_m)=T_{\sigma}\big(\varphi_R\ast f_1,\dots, \varphi_R\ast f_m \big).
\end{equation}
It follows from Theorem \ref{ann} with \eqref{conditionsnpn2} that
$$T_{\sigma_R}(f_1,\dots,f_m)\to T_{\sigma}(f_1,\dots,f_m) ~\text{ in the sense of tempered distribution}$$
as $R\to \infty$, because $\varphi_R$ is an approximation of identity (as $R\to \infty$).
This clearly deduces
$$
\lim_{R\to\infty} \phi_l \ast T_{\sigma_R}\big(f_1, \dots, f_m \big) (x) =  \phi_l \ast T_\sigma\big(f_1, \dots, f_m \big) (x) \q \text{ for } x\in \bbrn,
$$
and thus we have
\begin{equation}\label{newcondition2}
\begin{aligned}
\big\| T_\sigma (f_1, \dots, f_m) \big\|_{H^p(\bbrn)} &= \Big\| \sup_{l\in\Z} \big| \phi_l \ast T_\sigma (f_1, \dots, f_m) \big| \Big\|_{L^p(\bbrn)}\\
&\leq \Big\|  \liminf_{R\to\infty} \sup_{l\in\Z} \big| \phi_l \ast T_{\sigma_R} (f_1, \dots, f_m) \big|  \Big\|_{L^p(\bbrn)}\\
&\leq \liminf_{R\to\infty} \Big\| \sup_{l\in\Z} \big| \phi_l \ast T_{\sigma_R} (f_1, \dots, f_m) \big|  \Big\|_{L^p(\bbrn)}\\
&= \liminf_{R\to\infty}\big\Vert T_{\sigma_R}(f_1,\dots,f_m)\big\Vert_{H^p(\bbrn)}
\end{aligned}
\end{equation}
where Fatou's lemma is applied.
Moreover, using the compact support conditions of both $\wh{\Psi}$ and $\wh{\Phi_R}$, it is easy to see that
\begin{align*}
\LL_s^2[\sigma_R]&\le \sup_{j:2R<2^j\le 4R\sqrt{m}}\big\Vert \sigma(2^j\ccdot)\wh{\Psi}(\ccdot)\wh{\Phi_R}(2^j\ccdot)\big\Vert_{L^2_s((\bbrn)^m)}+\sup_{j:2^j\le 2R}\big\Vert \sigma(2^j\ccdot)\wh{\Psi}\big\Vert_{L_s^2((\bbrn)^m)}.
\end{align*}
The second term is clearly less than $\LL_s^2[\sigma]$ and the Kato-Ponce inequality ( see \cite{Gr_Oh2014} ) yields that for $2^j\sim R$
$$\big\Vert \sigma(2^j\ccdot)\wh{\Psi}(\ccdot)\wh{\Phi}(2^jR^{-1}\ccdot)\big\Vert_{L^2_s((\bbrn)^m)}\lesssim \big\Vert \sigma(2^j\ccdot)\wh{\Psi}\big\Vert_{L^2_s((\bbrn)^m)} +\big\Vert \sigma(2^j\ccdot)\wh{\Psi}\big\Vert_{L^2((\bbrn)^m)}\lesssim \LL_s^2[\sigma],$$
which finally proves
\begin{equation}\label{newcondition3}
\LL_s^2[\sigma_R]\lesssim \LL_s^2[\sigma] \q \text{ uniformly in } R>0.
\end{equation}

In view of \eqref{newcondition2} and \eqref{newcondition3}, we reduce the proof of Theorem \ref{main} to the estimate
\begin{equation}\label{mainreduction}
\big\Vert T_{\sigma_R}(f_1,\dots,f_m)\big\Vert_{H^p(\bbrn)}\lesssim \LL_s^2[\sigma_R]\prod_{\mu=1}^m\Vert f_{\mu}\Vert_{H^{p_{\mu}}(\bbrn)} \q \text{ uniformly in }~ R>0
\end{equation}
as this can be improved to any multipliers $\sigma$ in general:
\begin{equation}\label{generalapprox}
 \begin{aligned}
 \big\Vert T_{\sigma}(f_1,\dots,f_m)\big\Vert_{H^p(\bbrn)}&\le \liminf_{R\to \infty} \big\Vert T_{\sigma_R}(f_1,\dots,f_m)\big\Vert_{H^p(\bbrn)}\\
 &\lesssim \liminf_{R\to \infty}\LL_s^2[\sigma_R] \prod_{\mu=1}^m\Vert f_{\mu}\Vert_{H^{p_{\mu}}(\bbrn)}\\
 &\lesssim \LL_s^2[\sigma]\prod_{\mu=1}^m\Vert f_{\mu}\Vert_{H^{p_{\mu}}(\bbrn)}.
 \end{aligned}
\end{equation}

In order to establish \eqref{mainreduction}, we apply the atomic decomposition for Hardy spaces.
Each $f_{\mu}\in H^{p_{\mu}}(\bbrn)$ can be written as
$$f_{\mu}=\sum_{k_{\mu}\in\bbn}\la_{\mu,k_{\mu}}a_{\mu,k_{\mu}}$$
where $\{a_{\mu,k_{\mu}}\}_{k_{\mu}\in \bbn}$ is a sequence of $H^{p_{\mu}}$-atoms so that $\supp(a_{\mu,k_{\mu}})\subset Q_{\mu,k_{\mu}}$, $\Vert a_{\mu,k_{\mu}}\Vert_{L^{\infty}(\bbrn)}\le |Q_{\mu,k_{\mu}}|^{-1/p_{\mu}}$ for some cube $Q_{\mu,k_{\mu}}$, $\int_{Q_{\mu,k_{\mu}}} x^{\ga}a_{\mu,k_{\mu}}(x)dx=0$ for all multi-indices $|\ga|\le L$, and $\big( \sum_{k_{\mu}\in\bbn}|\la_{\mu,k_{\mu}}|^{p_{\mu}}\big)^{1/p_{\mu}}\lesssim \Vert f_{\mu}\Vert_{H^{p_{\mu}}(\bbrn)}$ where $L$ is a sufficiently large number greater than $n/p-n$.
Then we claim that
\begin{equation}\label{ptinequality}
\big| \phi_l\ast T_{\sigma_R}\big(f_1,\dots,f_m\big)(x)\big|\le \sum_{k_1,\dots,k_m\in\bbn}\Big(\prod_{\mu=1}^{m}|\la_{\mu,k_{\mu}}|\Big)\big| \phi_l\ast T_{\sigma_R}\big(a_{1,k_1},\dots,a_{m,k_m}\big)(x)\big|.
\end{equation}
For this one,  we define 
$$f_{\mu}^{N}:=\sum_{k_{\mu}=1}^{N}\la_{\mu,k_{\mu}}a_{\mu,k_{\mu}}, \qq N\in \bbn$$
so that 
\begin{equation}\label{hpapproximate}
\big\Vert f_{\mu}-f_{\mu}^{N}\big\Vert_{H^{p_{\mu}}(\bbrn)}\lesssim \Big( \sum_{k_{\mu}=N+1}^{\infty}\big| \la_{\mu,k_{\mu}}\big|^{p_{\mu}}\Big)^{1/p_{\mu}}\to 0 \q \text{ as }~ N\to \infty.
\end{equation}
Then the left-hand side of \eqref{ptinequality} is less than
\begin{align*}
&\big| \phi_l\ast T_{\sigma_R}\big(f_1-f_1^{N},f_2,\dots,f_m\big)(x)\big| \\
&~ + \big| \phi_l\ast T_{\sigma_R}\big(f_1^{N},f_2-f_2^{N},f_3,\dots,f_m\big)(x)\big|\\
&\qq \qq \vdots \\
&\qq + \big| \phi_l\ast T_{\sigma_R}\big(f_1^{N},\dots,f_{m-1}^{N},f_m-f_m^{N}\big)(x)\big|\\
&\qq \q +\big| \phi_l\ast T_{\sigma_R}\big(f_1^{N},\dots,f_m^{N}\big)(x)\big|.
\end{align*}
The last term is clearly dominated by
\begin{align*}
&\sum_{k_1=1}^{N}\cdots\sum_{k_{m}=1}^{N}\Big(\prod_{\mu=1}^{m}|\la_{\mu,k_{\mu}}|\Big)\big| \phi_l\ast T_{\sigma_R}\big(a_{1,k_1},\dots,a_{m,k_m}\big)(x)\big|\\
&\le \sum_{k_1,\dots,k_m\in\bbn}\Big(\prod_{\mu=1}^{m}|\la_{\mu,k_{\mu}}|\Big)\big| \phi_l\ast T_{\sigma_R}\big(a_{1,k_1},\dots,a_{m,k_m}\big)(x)\big|
\end{align*} by the multi-linearity of the operator $\phi_l\ast T_{\sigma_R}$. 
Therefore, it suffices to show that the remaining terms approach $0$ as $N\to \infty$. We will be concerned only with the first term as similar arguments can be applied to the others.
We note that the Fourier transform of $T_{\sigma_R}(f_1-f_1^N,f_2,\dots,f_m)$ is supported on the ball $\{\xi\in\bbrn: |\xi|\le 2\sqrt{m}R\}$ due to the support of $\wh{\varphi_R}$. This shows that
$$ \phi_l\ast T_{\sigma_R}\big(f_1-f_1^{N},f_2,\dots,f_m\big)=\phi_l\ast \Upsilon_R\ast T_{\sigma_R}\big(f_1-f_1^{N},f_2,\dots,f_m\big)$$
for some Schwartz function $\Upsilon_R$ on $\bbrn$ whose Fourier transform is equal to $1$ on the ball $\{\xi\in\bbrn: |\xi|\le 2\sqrt{m}R\}$ and is supported on the larger ball
$\{\xi\in\bbrn: |\xi|\le 4\sqrt{m}R\}$. Now, using Nikolskii's inequality (sometimes called Bernstein's inequality) in \cite[Proposition 1.3.2]{Tr1983},
we have
\begin{align*}
&\big| \phi_l\ast T_{\sigma_R}\big(f_1-f_1^{N},f_2,\dots,f_m\big)(x)\big|\\
&=\Big| \int_{\bbrn} \phi_l\ast \Upsilon_R(x-y)T_{\sigma_R}\big(f_1-f_1^N,f_2,\dots,f_m \big)(y)  dy\Big|\\
&\lesssim R^{n/p-n} \Big( \int_{\bbrn} \big| \phi_l\ast \Upsilon_R(x-y)\big|^p \big|T_{\sigma_R}\big(f_1-f_1^N,f_2,\dots,f_m \big)(y) \big|^p dy\Big)^{1/p}\\
&\lesssim_{R,l,p}\big\Vert T_{\sigma_R}(f_1-f_1^N,f_2,\dots,f_m ) \big\Vert_{L^p(\bbrn)}\\
&\lesssim \LL_s^2[\sigma_R]\big\Vert f_1-f_1^N\big\Vert_{H^{p_1}(\bbrn)}\prod_{\mu=2}^{m}\Vert f_{\mu}\Vert_{H^{p_{\mu}}(\bbrn)} \to 0
\end{align*}
as $N\to \infty$, thanks to \eqref{hpapproximate}. Here, the last inequality follows from Theorem \ref{ann}. This completes the proof of \eqref{ptinequality}.

We also remark that the vanishing moment condition \eqref{multivanishing}, together with \eqref{sigmarrepresentation}, guarantees
\begin{equation}\label{newcondition1}
\int_{\bbrn}x^{\alpha}T_{\sigma_R}\big(a_{1,k_1},\dots,a_{m,k_m}\big)(x)dx=0 \q \text{ for multi-indices }~ |\alpha|\le \frac{n}{p}-n
\end{equation} as each $\varphi_R\ast a_{\mu,k_{\mu}}$ belongs to $\mathscr{S}^L(\bbrn)$.

In conclusion, matters reduce to the following proposition, in view of \eqref{hardydef2}, \eqref{ptinequality}, and \eqref{newcondition1}.
\begin{proposition}\label{mainpropo}
Let $0<p_1,\dots,p_m\le 1$ and $0<p\le 1/m$ with $1/p=1/p_1+\dots+1/p_m$.
Suppose that $L$ is sufficiently large and
$$s>\frac{n}{p}-\frac{n}{2}.$$
We also assume that the multilinear multiplier operator $T_{\sigma}$ satisfies the vanishing moment condition
\begin{equation*}
\int_{\bbrn}x^{\alpha}T_{\sigma}\big(a_{1,k_1},\dots,a_{m,k_m}\big)(x) dx=0 \q \text{ for multi-indices }~ |\alpha|\le \frac{n}{p}-n
\end{equation*}
for functions $f_{\mu}\in \mathscr{S}^L(\bbrn)$ having the atomic representation $f_{\mu}=\sum_{k_{\mu}\in\bbn}\la_{\mu,k_{\mu}}a_{\mu,k_{\mu}}$, $\mu=1,\dots,m$, where each $a_{\mu,k_{\mu}}$ is an $H^{p_{\mu}}(\bbrn)$ atom with the property \eqref{hardyatomvanishing} upto the degree $L$, and $\{\la_{\mu,k_{\mu}}\}_{k_{\mu}\in \bbn}$ is a sequence of coefficients satisfying $\big( \sum_{k_m\in\bbn}|\la_{\mu,k_{\mu}}|^{p_{\mu}}\big)^{1/p{\mu}}\lesssim \Vert f_{\mu}\Vert_{H^{p_{\mu}}(\bbrn)}$.
Then we have
\begin{equation*}
\begin{aligned}
&\bigg\Vert \sup_{l\in\bbz}\sum_{k_1,\dots,k_m\in\bbn}\Big(\prod_{\mu=1}^{m}|\la_{\mu,k_{\mu}}|\Big)\big| \phi_l\ast T_{\sigma}(a_{1,k_1},\dots,a_{m,k_m})\big|     \bigg\Vert_{L^p(\bbrn)}\lesssim \LL_s^2[\sigma]\prod_{\mu=1}^{m}\Vert f_{\mu}\Vert_{H^{p_{\mu}}(\bbrn)}.
\end{aligned}
\end{equation*}
\end{proposition}
The proof of the proposition will be given in the next section.

\section{Proof of Theorem \ref{main} : Proof of Proposition \ref{mainpropo}}

Suppose that each $H^{p_{\mu}}$ atom $a_{\mu,k_{\mu}}$ is associated with a cube $Q_{\mu,k_{\mu}}$ so that $\supp(a_{\mu,k_{\mu}})\subset Q_{\mu,k_{\mu}}$, $\Vert a_{\mu,k_{\mu}}\Vert_{L^{\infty}(\bbrn)}\le |Q_{\mu,k_{\mu}}|^{-1/p_{\mu}}$, and $\int_{Q_{\mu,k_{\mu}}}x^{\ga}a_{\mu,k_{\mu}}(x)dx=0$ for multi-indices $|\ga|\le L$.
Then we write 
$$\sum_{k_1,\dots,k_m\in\bbn}\Big(\prod_{\mu=1}^{m}|\la_{\mu,k_{\mu}}|\Big)\big| \phi_l\ast T_{\sigma}\big(a_{1,k_1},\dots,a_{m,k_m}\big)(x)\big| \le \mathcal{G}_{\mathrm{in}}^l(x)+\mathcal{G}_{\mathrm{out}}^l(x)$$  
where
$$\mathcal{G}_{\mathrm{in}}^l(x):=\sum_{k_1,\dots,k_m\in\bbn}\Big( \prod_{\mu=1}^{m}|\la_{\mu,k_{\mu}}|\Big)\big| \phi_l\ast T_{\si}\big(a_{1,k_1},\dots,a_{m,k_m}\big)(x) \big|\chi_{Q_{1,k_1}^{**}\cap\cdots\cap Q_{m,k_m}^{**}}(x),$$
$$\mathcal{G}_{\mathrm{out}}^l(x):=\sum_{k_1,\dots,k_m\in\bbn}\Big( \prod_{\mu=1}^{m}|\la_{\mu,k_{\mu}}|\Big)\big|\phi_l\ast T_{\si}\big(a_{1,k_1},\dots,a_{m,k_m}\big)(x) \big|\chi_{(Q_{1,k_1}^{**}\cap\cdots\cap Q_{m,k_m}^{**})^c}(x).$$
Now it suffices to show that
\begin{equation*}
\Big\Vert \sup_{l\in\bbz}\mathcal{G}_{\mathrm{in}/\mathrm{out}}^l\Big\Vert_{L^p(\bbrn)} \lesssim \LL_s^2[\sigma] \prod_{\mu=1}^{m}\Vert f_{\mu}\Vert_{H^{p_{\mu}}(\bbrn)}.
\end{equation*}

\hfill

The term associated with $\mathcal{G}_{\mathrm{in}}^l$ can be easily handled by the argument of Grafakos and Kalton \cite{Gr_Ka2001}. 
Suppose that $Q_{1,k_1}^{**}\cap \cdots\cap Q_{m,k_m}^{**}\not=\emptyset$ as we are done if the intersection is empty. We may also assume that the side length of the cube $Q_{1,k_1}$ is the smallest among those of the cubes $Q_{1,k_1},\dots,Q_{m,k_m}$. Then as done in \cite{Gr_Ka2001}, we choose a cube $R_{\kkk}$, $\kkk:=(k_1,\dots,k_m)$, such that 
\begin{equation*}
Q_{1,k_1}^{**} \cap \cdots \cap Q_{m,k_m}^{**} \subset R_{\kkk}  \subset Q_{1,k_1}^{***} \cap \cdots \cap Q_{m,k_m}^{***}
\end{equation*}
and $\ell(R_{\kkk})\gtrsim \ell(Q_{1, k_1})$.
By applying Lemma \ref{reverse FS ineq}, we have
\begin{align*}
\Big\Vert \sup_{l\in\bbz} \mathcal{G}_{\mathrm{in}}^l \Big\Vert_{L^p(\bbrn)}&\lesssim 
\bigg\|   \sum_{k_1, \cdots, k_m\in \bbn} 
	\Big( \frac{1}{ |R_{\kkk}| } \int_{R_{\kkk}} \big| \mathcal{M}T_\sigma\big(a_{1,k_1},\cdots, a_{m,k_m}\big)(y)\big| dy\Big)\\
&\qq\qq\qq\qq\qq\qq\times \Big( \prod_{\mu=1}^m |\lmd_{\mu,k_{\mu}}| \Big)	\chi_{R_{\kkk}} \bigg\|_{L^p(\bbrn)}
\end{align*}
and the integral over $R_{\kkk}$ is less than
\begin{align*}
&|R_{\kkk}|^{1/2} \big\Vert \mathcal{M}T_\sigma\big(a_{1,k_1},\cdots, a_{m,k_m}\big)\big\Vert_{L^2(\bbrn)}\\
&\lesssim |R_{\kkk}|^{1/2} \LL_{s}^2[\sigma]\Vert a_{1,k_1}\Vert_{L^2(\bbrn)}\prod_{\mu=2}^{m}\Vert a_{\mu,k_{\mu}}\Vert_{L^{\infty}(\bbrn)}\\
&\le |R_{\kkk}|^{1/2} \LL_{s}^2[\sigma]|Q_{1,k_1}|^{1/2}\prod_{\mu=1}^{m}|Q_{\mu,k_{\mu}}|^{-1/p_{\mu}}\lesssim  |R_{\kkk}|\LL_{s}^2[\sigma]\prod_{\mu=1}^{m}|Q_{\mu,k_{\mu}}|^{-1/p_{\mu}}
\end{align*}
where we applied the maximal inequality for $\mathcal{M}$ and the $L^2\times L^{\infty}\times \cdots \times L^{\infty}$ boundedness of $T_{\sigma}$, which holds due to Theorem \ref{ann}, in the first inequality.
This, together with H\"older's inequality, finally concludes that
\begin{align*}
\Big\Vert \sup_{l\in\bbz}  \mathcal{G}_{\mathrm{in}}^l \Big\Vert_{L^p(\bbrn)}&\lesssim \LL_s^2[\sigma]\bigg\Vert \prod_{\mu=1}^{m}\Big( \sum_{k_{\mu}\in\bbn}|\la_{\mu,k_{\mu}}||Q_{\mu,k_{\mu}}|^{-1/p_{\mu}}\chi_{Q_{\mu,k_{\mu}}^{***}}\Big)\bigg\Vert_{L^p(\bbrn)}\\
&\lesssim \LL_s^2[\sigma]\prod_{\mu=1}^{m}\Big\Vert \sum_{k_{\mu}\in\bbn}|\la_{\mu,k_{\mu}}||Q_{\mu,k_{\mu}}|^{-1/p_{\mu}}\chi_{Q_{\mu,k_{\mu}}^{***}}\Big\Vert_{L^{p_{\mu}}(\bbrn)}\\
&\lesssim \LL_s^2[\sigma]\prod_{\mu=1}^{m}\Big( \sum_{k_{\mu}\in\bbn}\big| \la_{\mu,k_{\mu}}\big|^{p_{\mu}}\Big)^{1/p_{\mu}}         \lesssim \LL_s^2[\sigma]\prod_{\mu=1}^{m}\Vert f_{\mu}\Vert_{H^{p_{\mu}}(\bbrn)}
\end{align*} as desired.\\

To estimate the remaining term associated with $\mathcal{G}_{\mathrm{out}}^l$, we claim the following lemma.
\begin{lemma}\label{keylemma}
Let $0<p_1,\dots,p_m\le 1$ and $0<p\le 1/m$ with $1/p_1+\dots+1/p_m=1/p$. 
Suppose that $L$ is a sufficiently large number and let $a_{\mu}$, $1\le \mu\le m$, be $H^{p_{\mu}}$-atoms such that 
$$\mathrm{supp}(a_{\mu})\subset Q_{\mu},\q \Vert a_{\mu}\Vert_{L^{\infty}(\bbrn)}\le |Q_{\mu}|^{-1/p_{\mu}},\q \int_{Q_{\mu}}x^{\alpha}a_{\mu}(x)dx=0$$
for all $|\alpha|\le L$.
For each $\mu=1,\dots,m$, let
$$B_{\mu,l}:=\big\{x\in\bbrn: |x-c_{Q_{\mu}}|\le 100n^22^{-l} \big\}$$
where $c_{Q_{\mu}}$ is the center of the cube $Q_{\mu}$.
Let $J_0$ be any nonempty subset of $J_m:=\{1,\dots,m\}$ and let $I_0$ be any nonempty subset of $J_0$. 
Suppose that $s>n/p-n/2$ and $\sigma$ satisfies the vanishing moment condition
\begin{equation}\label{atomicvanishing}
\int_{\bbrn}x^{\alpha}T_{\sigma}\big(a_1,\dots,a_m\big)(x) dx=0 \q \text{ for multi-indices }~ |\alpha|\le \frac{n}{p}-n.
\end{equation} 
 Then the following two statements are valid:
\begin{enumerate}
\item There exist nonnegative functions $\Ga_1^{J_0},\dots,\Ga_m^{J_0}$ such that
$$\Vert \Ga_{\mu}^{J_0}\Vert_{L^{p_{\mu}}(\bbrn)}\lesssim 1 \q \text{ for all }~1\le \mu\le m,$$ 
and
\begin{equation}\label{deresult1}
\big| \phi_l\ast T_{\sigma}\big(a_1,\dots,a_m\big)(x)\big|\lesssim \LL_s^2[\sigma]\Ga_1^{J_0}(x)\cdots \Ga_m^{J_0}(x)\q \text{ uniformly in }~ l\in\bbz
\end{equation}
for all $x\in \Big[ \big( \bigcap_{\mu \in J_m\setminus J_0}Q_{\mu}^{**}\big)\setminus \big(\bigcup_{\mu\in J_0}Q_{\mu}^{**}\big)\Big]\cap \big( \bigcap_{\mu\in J_0}B_{\mu,l}\big)$.

\item There exist nonnegative functions $\La_1^{J_0},\dots,\La_m^{J_0}$ such that
$$\Vert \La_{\mu}^{J_0}\Vert_{L^{p_{\mu}}(\bbrn)}\lesssim 1 \q \text{ for all }~1\le \mu\le m,$$ 
and
\begin{equation}\label{keylemma2assertion}
\big| \phi_l\ast T_{\sigma}\big(a_1,\dots,a_m\big)(x)\big|\lesssim \LL_s^2[\sigma]\La_1^{J_0}(x)\cdots \La_m^{J_0}(x)\q \text{ uniformly in }~ l\in\bbz
\end{equation}
for all $x\in \Big[ \big( \bigcap_{\mu\in J_m\setminus J_0}Q_{\mu}^{**}\big)\setminus \big(\bigcup_{\mu \in J_0}Q_{\mu}^{**}\big)\Big]\cap \Big[ \big( \bigcap_{\mu\in J_0\setminus I_0}B_{\mu,l}\big)\setminus \big(\bigcup_{\mu\in I_0}B_{\mu,l}\big)\Big]$.
\end{enumerate}
Here, we adopt the convention that an intersection over an empty set is considered $\bbrn$.
\end{lemma}
The proof of the lemma will be provided in the next section.
Taking the lemma for granted and using the fact that
\begin{align*}
\Big( \bigcap_{\mu=1}^{m}Q_{\mu,k_{\mu}}^{**}\Big)^c& =\bigcup_{J_0\subset J_m, J_0\not=\emptyset} \Big[  \Big( \bigcap_{\mu \in J_m\setminus J_0}Q_{\mu}^{**} \Big) \cap \Big(\bigcap_{\mu\in J_0} \big( Q_{\mu}^{**} \big)^c \Big) \Big] \\
&=\bigcup_{J_0\subset J_m, J_0\not=\emptyset} \Big[ \Big( \bigcap_{\mu \in J_m\setminus J_0}Q_{\mu}^{**}\Big)\setminus \Big(\bigcup_{\mu\in J_0}Q_{\mu}^{**}\Big) \Big],
\end{align*}
we have
$$\big| \phi_l\ast T_{\sigma}\big(a_{1,k_{1}},\dots,a_{m,k_m}\big)(x)\big|\chi_{(\bigcap_{\mu=1}^{m}Q_{\mu,k_{\mu}}^{**})^c}(x)\lesssim \LL_s^2[\sigma]\prod_{\mu=1}^{m}\Big[\sum_{J_0\subset J_m, J_0\not= \emptyset}\big(\Ga_{\mu}^{J_0}(x)+ \La_{\mu}^{J_0}(x) \big)\Big]$$
and this implies that
\begin{align*}
\Big\Vert \sup_{l\in\bbz} \mathcal{G}_{\mathrm{out}}^l \bigg\Vert_{L^p(\bbrn)}&\lesssim \LL_s^2[\sigma]\bigg\Vert\prod_{\mu=1}^{m}\Big( \sum_{k_{\mu}\in\bbn} |\la_{\mu,k_{\mu}}|\sum_{J_0\subset J_m, J_0\not= \emptyset}\big(\Ga_{\mu}^{J_0}+ \La_{\mu}^{J_0} \big)\Big)\bigg\Vert_{L^p(\bbrn)}\\
&\le \LL_s^2[\sigma]\prod_{\mu=1}^{m}\bigg\Vert \sum_{k_{\mu}\in\bbn} |\la_{\mu,k_{\mu}}| \sum_{J_0\subset J_m, J_0\not= \emptyset}\big(\Ga_{\mu}^{J_0}+ \La_{\mu}^{J_0} \big)\bigg\Vert_{L^{p_{\mu}}(\bbrn)}
\end{align*}
by using H\"older's inequality.
In addition,  it follows that for each $\mu=1,\dots,m$,
\begin{align*}
&\bigg\Vert \sum_{k_{\mu}\in\bbn} |\la_{\mu,k_{\mu}}|\sum_{J_0\subset J_m, J_0\not= \emptyset}\big(\Ga_{\mu}^{J_0}+ \La_{\mu}^{J_0} \big)\bigg\Vert_{L^{p_{\mu}}(\bbrn)}\\
&\le \bigg( \sum_{k_{\mu}\in\bbn} |\la_{\mu,k_{\mu}}|^{p_{\mu}} \sum_{J_0\subset J_m, J_0\not= \emptyset}  \Big(  \big\Vert\Ga_{\mu}^{J_0} \big\Vert_{L^{p_{\mu}}(\bbrn)}^{p_{\mu}}+\big\Vert\La_{\mu}^{J_0} \big\Vert_{L^{p_{\mu}}(\bbrn)}^{p_{\mu}}\Big)    \bigg)^{1/p_{\mu}}\\
&\lesssim \Big( \sum_{k_{\mu}\in\bbn}|\la_{\mu,k_{\mu}}|^{p_{\mu}}\Big)^{1/p_{\mu}}\lesssim \Vert f_{\mu}\Vert_{H^{p_{\mu}}(\bbrn)}
\end{align*}
and this deduces
$$\Big\Vert \sup_{l\in\bbz} \mathcal{G}_{\mathrm{out}}^l \bigg\Vert_{L^p(\bbrn)}\lesssim \LL_s^2[\sigma]\prod_{\mu=1}^{m}\Vert f_{\mu}\Vert_{H^{p_{\mu}}(\bbrn)}.$$

This completes the proof of Proposition \ref{mainpropo}.

\section{Proof of Lemma \ref{keylemma}}\label{keylemmasection}

Without loss of generality, we may assume that $J_0=\{1,\dots,\nu\}$ for some $1\le \nu\le m$.
 Let $[r]$ denote the greatest integer less than or equal to $r \in \mathbb{R}$.
Since $s>n/p-n/2$, there exists $0<\delta<1$ such that 
$$n/p-n<[n/p-n]+\delta<s-n/2,$$
 and then we choose $0<\epsilon<1/(m+1)$ satisfying
 \begin{equation}\label{conditionons}
 n/p-n+(m+1)\epsilon<[n/p-n]+\delta<s-n/2,
 \end{equation}
 which implies that
 \begin{equation}\label{npnme}
 [n/p-n+(m+1)\epsilon]\le [n/p-n]\le n/p-n
 \end{equation}
 because $0<\delta<1$.
We may also assume that 
\begin{equation}\label{keycondition2}
\ell(Q_1)=\min_{1\le \ii\le \nu}\{\ell(Q_{\ii}) \}.
\end{equation}

By using the Littlewood-Paley decomposition with $\{ \Psi_j\}_{j\in\bbz}$, 
we can write
\begin{equation}\label{lpdecomp}
T_{\sigma}(a_1,\dots,a_m)=\sum_{j\in\bbz}T_j(a_1,\dots,a_m)
\end{equation}
where
$$K_j(\yyy):=\big( \sigma(2^j\ccdot)\wh{\Psi}\big)^{\vee}(\yyy)$$
and
$$T_j\big(a_1,\dots,a_m\big)(x):=2^{jmn}\int_{(\bbrn)^m}  K_j\big(2^j(x-y_1),\dots,2^j(x-y_m) \big)f_1(y_1)\cdots f_m(y_m)\;    d\yyy .$$

\subsection{The first assertion}
Suppose that 
\begin{equation}\label{rangeofx}
x\in  \Big( \bigcap_{\mu =1}^{\nu}(Q_{\mu}^{**})^c\cap B_{\mu,l}\Big)\medcap \Big(\bigcap_{\mu\in J_m\setminus J_0}Q_{\mu}^{**}\Big).
\end{equation}
We will prove \eqref{deresult1} by setting
$$\Ga_{\mu}^{J_0}(x):=\begin{cases}
\frac{\ell(Q_{\mu})^{\epsilon}}{|x-c_{Q_{\mu}}|^{n/p_{\mu}+\epsilon}}\chi_{(Q_{\mu}^{**})^c}(x) & \mu\in J_0\\
\ell(Q_{\mu})^{-n/p_{\mu}}\chi_{Q_{\mu}^{**}}(x) & \mu\in J_m\setminus J_0
\end{cases}.$$
Then it is clear that $\big\Vert \Ga_{\mu}^{J_0}\big\Vert_{L^{p_{\mu}}(\bbrn)}\lesssim 1$ for all $\mu\in J_m$.
In order to establish \eqref{deresult1},
let $$N_{\nu}:=[n/p_1+\dots+n/p_{\nu}+\nu\epsilon-n]$$ and $$\epsilon_{\nu}:=n/p_1+\dots+n/p_{\nu}+\nu\epsilon-n-N_{\nu}$$
so that
\begin{equation}\label{keycondition1}
N_{\nu}+n+\epsilon_{\nu}=\sum_{\mu=1}^{\nu}\big(n/p_{\mu}+\epsilon\big).
\end{equation}
We note that \eqref{npnme} yields 
\begin{equation}\label{conditionn}
N_{\nu}\le N_{m}\le n/p-n,
\end{equation}
and  for any $1\le \mu\le \nu$,
$$2^{l}\lesssim \frac{1}{|x-c_{Q_{\mu}}|}$$
because of \eqref{rangeofx}.
By using the vanishing moment condition \eqref{atomicvanishing}, together with \eqref{conditionn} and Lemma \ref{eps control}, we write 
\begin{align*}
\big| \phi_l\ast T_{\sigma}\big(a_1,\dots,a_m\big)(x)\big| &\lesssim 2^{l(N_{\nu}+n+\epsilon_{\nu})}\int_{\bbrn}|y-c_{Q_1}|^{N_{\nu}+\epsilon_{\nu}}\big|T_{\sigma}\big(a_1,\dots,a_m\big)(y) \big| \;dy\\
&\lesssim \Big( \prod_{\mu=1}^{\nu}\frac{1}{|x-c_{Q_{\mu}}|^{n/p_{\mu}+\epsilon}}\Big)\int_{\bbrn} |y-c_{Q_1}|^{N_{\nu}+\epsilon_{\nu}}\big| T_{\sigma}\big(a_1,\dots,a_m\big)(y)\big| \; dy\\
&\le \II_1(x)+\II_2(x)
\end{align*}
where
$$\II_1(x):=\Big( \prod_{\mu=1}^{\nu}\frac{\chi_{(Q_{\mu}^{**})^c}(x)}{|x-c_{Q_{\mu}}|^{n/p_{\mu}+\epsilon}}\Big)\Big(\prod_{\mu\in J_m\setminus J_0}\chi_{Q_{\mu}^{**}}(x) \Big) \int_{Q_1^*} |y-c_{Q_1}|^{N_{\nu}+\epsilon_{\nu}}\big| T_{\sigma}\big(a_1,\dots,a_m\big)(y)\big| \; dy$$
and
$$\II_2(x):=  \Big( \prod_{\mu=1}^{\nu}\frac{\chi_{(Q_{\mu}^{**})^c}(x)}{|x-c_{Q_{\mu}}|^{n/p_{\mu}+\epsilon}}\Big)\Big(\prod_{\mu\in J_m\setminus J_0}\chi_{Q_{\mu}^{**}}(x) \Big) \int_{(Q_1^*)^c}|y-c_{Q_1}|^{N_{\nu}+\epsilon_{\nu}}\big| T_{\sigma}\big(a_1,\dots,a_m\big)(y)\big| \; dy.$$
Now it suffices to show that
\begin{equation}\label{keyclaim}
\II_{\ii}(x)\lesssim \LL_s^2[\sigma]\Ga^{J_0}_1(x)\dots\Ga^{J_0}_m(x), \qq \ii=1,2.
\end{equation}

To estimate $\II_1$, we see that
\begin{align*}
& \int_{Q_1^*} |y-c_{Q_1}|^{N_{\nu}+\epsilon_{\nu}}\big| T_{\sigma}\big(a_1,\dots,a_m\big)(y)\big| \; dy \lesssim \ell(Q_1)^{N_{\nu}+\epsilon_{\nu}}\big\Vert T_{\sigma}(a_1,\dots,a_m)\big\Vert_{L^1(\bbrn)}\\
 &\lesssim\LL_s^2[\sigma]\ell(Q_1)^{N_{\nu}+\epsilon_{\nu}}\Vert a_1\Vert_{H^1(\bbrn)}\prod_{\mu=2}^{m}\Vert a_{\mu}\Vert_{L^{\infty}(\bbrn)}   \le \LL_s^2[\sigma]\ell(Q_1)^{N_{\nu}+n+\epsilon_{\nu}}\prod_{\mu=1}^{m}\ell(Q_{\mu})^{-n/p_{\mu}}     \\
& \lesssim \LL_s^2[\sigma]\Big(\prod_{\mu=1}^{\nu}\ell(Q_{\mu})^{\epsilon} \Big) \Big( \prod_{\mu\in J_m\setminus J_0}\ell(Q_{\mu})^{-n/p_{\mu}}\Big)
\end{align*}
where we applied the $H^1\times L^{\infty}\times \cdots\times L^{\infty}\to L^1$ boundedness of $T_{\sigma}$ in Theorem \ref{ann} and the conditions \eqref{keycondition1}, \eqref{keycondition2}.
Here, we notice that our assumption $s>\frac{n}{p} - \frac{n}{2}\geq mn - \frac{n}{2}$ implies $s > \frac{mn}{2}$ whenever $m\geq1$, and  $|Q_1|^{1/p_1-1}a_1$ is an $H^1$-atom for which $\Vert a_1\Vert_{H^1(\bbrn)}\lesssim |Q_1|^{1-1/p_1}$.
This proves \eqref{keyclaim} for $\ii=1$.\\

For the other term $\II_2$, we decompose $T_{\sigma}$ as in (\ref{lpdecomp}) and then investigate the behavior of the operators $T_j$.
We actually claim that for $y\in (Q_1^*)^c$
\begin{align}\label{keyestimate1}
\big| T_j\big(a_1,\dots,a_m\big)(y)\big|&\lesssim 2^{j\nu n} \Big(\prod_{\mu=1}^{\nu}\ell(Q_{\mu})^{-n/p_{\mu}+n}\Big)\Big( \prod_{\mu\in J_m\setminus J_0}\ell(Q_{\mu})^{-n/p_{\mu}}\Big)\\
&\qq\qq\times \frac{1}{\langle 2^j(y-c_{Q_1})\rangle^{N_{\nu}+\epsilon_{\nu}+n/2+\epsilon}}\min\big\{ h_{j,s}^{Q_1,0}(y),h_{j,s}^{Q_1,L}(y)\big\} \nonumber
\end{align}
where
\begin{align}\label{hjs1def}
h_{j,s}^{Q_1,0}(y)&:=\int_{Q_1} |Q_1|^{-1} \Big\Vert \big\langle 2^j(y-z_1),z_2,\dots,z_m \big\rangle^{s} \big|K_j\big(2^j(y-z_1),z_2,\dots,z_m\big) \big|\Big\Vert_{L^2(z_2,\dots,z_m)} dz_1
\end{align}
and 
\begin{align}\label{hjs1ldef}
h_{j,s}^{Q_1,L}(y)&:=\big( 2^j\ell(Q_1)\big)^L\sum_{|\alpha|=L}\int_0^1 \bigg( \int_{Q_1} |Q_1|^{-1}   \\
&\qq\times \Big\Vert \big\langle 2^jy_{c_{Q_1},z_1}^{\theta},z_2,\dots,z_m \big\rangle^{s} \big|K_j\big(2^jy_{c_{Q_1},z_1}^{\theta},z_2,\dots,z_m\big) \big|\Big\Vert_{L^2(z_2,\dots,z_m)} dz_1 \bigg)d\theta \nonumber
\end{align}
by setting $y_{c_{Q_1},z_1}^{\theta}:=y-c_{Q_1}-\theta(z_1-c_{Q_1})$.
Here, it is easy to see by Minkowski's inequality that
\begin{equation}\label{hjl2norms}
\big\Vert h_{j,s}^{Q_1,0} \big\Vert_{L^2(\bbrn)} \lesssim 2^{-jn/2}\LL_{s}^2[\sigma], \qq \big\Vert h_{j,s}^{Q_1,L} \big\Vert_{L^2(\bbrn)} \lesssim 2^{-jn/2}\big(2^j\ell(Q_1) \big)^L\LL_{s}^2[\sigma].
\end{equation}  
Once the claim (\ref{keyestimate1}) is established, it follows that
\begin{align*}
&\int_{(Q_1^{*})^c}  |y-c_{Q_1}|^{N_{\nu}+\epsilon_{\nu}}\big| T_j\big(a_1,\dots,a_m\big)(y)\big|    dy\nonumber\\
&\lesssim 2^{-j(N_{\nu}+\epsilon_{\nu}-\nu n)}\Big(\prod_{\mu=1}^{\nu}\ell(Q_{\mu})^{-n/p_{\mu}+n} \Big)\Big(\prod_{\mu\in J_m\setminus J_0}\ell(Q_{\mu})^{-n/p_{\mu}} \Big)\\
&\qq\qq\qq\times \int_{(Q_1^*)^c} \frac{1}{\langle 2^j(y-c_{Q_1})\rangle^{n/2+\epsilon}}\min \big\{h_{j,\nu}^{Q_1,0}(y),h_{j,\nu}^{Q_1,L}(y) \big\} dy \nonumber
\end{align*}
and the integral in the last displayed expression is, via the Cauchy-Schwarz inequality, less than
\begin{align*}
&\bigg\Vert \frac{1}{\langle 2^j(\cdot-c_{Q_1})\rangle^{n/2+\epsilon}}\bigg\Vert_{L^2((Q_1^*)^c)}\min\Big\{\big\Vert h_{j,\nu}^{Q_1,0}\big\Vert_{L^2(\bbrn)}, \big\Vert h_{j,\nu}^{Q_1,L}\big\Vert_{L^2(\bbrn)}\Big\}\\
&\lesssim \LL_s^2[\sigma]\ell(Q_1)^{-\epsilon}2^{-jn-j\epsilon}\min\big\{1,\big(2^j\ell(Q_1)\big)^L \big\}
\end{align*}
by applying \eqref{hjl2norms}. Therefore, in view of \eqref{lpdecomp}, we have
\begin{align*}
&\int_{(Q_1^{*})^c}  |y-c_{Q_1}|^{N_{\nu}+\epsilon_{\nu}}\big| T_{\sigma}\big(a_1,\dots,a_m\big)(y)\big|    dy\\
&\lesssim \LL_s^2[\sigma]\Big(\prod_{\mu=1}^{\nu}\ell(Q_{\mu})^{-n/p_{\mu}+n} \Big)\Big(\prod_{\mu\in J_m\setminus J_0}\ell(Q_{\mu})^{-n/p_{\mu}} \Big)\\
&\qq\qq\qq\times \ell(Q_1)^{-\epsilon}\sum_{j\in\bbz}2^{-j\big(\epsilon+\sum_{\mu=1}^{\nu}(n/p_{\mu}-n+\epsilon) \big)}\min \Big\{1,\big(2^j\ell(Q_1)\big)^L \Big\}\\
&\lesssim \LL_s^2[\sigma]\Big(\prod_{\mu=1}^{\nu}\ell(Q_{\mu})^{\epsilon} \Big)\Big(\prod_{\mu\in J_m\setminus J_0}\ell(Q_{\mu})^{-n/p_{\mu}} \Big)
\end{align*}
due to \eqref{keycondition1}, \eqref{keycondition2}, and the fact that $$L>n/p-n\ge \epsilon+\sum_{\mu=1}^{\nu}(n/p_{\mu}-n+\epsilon) \q\text{ for all } \nu=1,\dots,m.$$
This finally proves \eqref{keyclaim} for $\ii=2$.

\hfill

Therefore, it remains to show the estimate \eqref{keyestimate1}. For notational convenience, let $$M_{\nu}^{\epsilon}:=N_{\nu}+\epsilon_{\nu}+n/2+\epsilon.$$
Moreover, for any $r>0$ and multi-indices $\alpha\in (\bbn_0)^n$,
we simply write
\begin{equation}\label{kjrdef}
\begin{aligned}
K_j^r(z_1,\dots,z_m)&:=\langle z_1,\dots,z_m \rangle^{r}K_j(z_1,\dots,z_m),\\
\big(\partial_1^{\alpha}K_j\big)^r(z_1,\dots,z_m)&:=\langle z_1,\dots,z_m \rangle^{r}\partial_1^{\alpha}K_j(z_1,\dots,z_m).
\end{aligned}
\end{equation}
We now observe that for $y\in (Q_1^*)^c$ and $z_1\in Q_1$, 
\begin{equation}\label{ycqyz}
|y-c_{Q_1}|\lesssim |y-z_1|
\end{equation}
and this implies that
\begin{align*}
&\langle 2^j(y-c_{Q_1})\rangle^{M_{\nu}^{\epsilon}}\big| T_j\big(a_1,\dots,a_m\big)(y)\big|\\
&\lesssim 2^{jmn}\Big(\prod_{\mu=1}^{m}\ell(Q_{\mu})^{-n/p_{\mu}}\Big)\int_{Q_1\times\cdots\times Q_m}\big|K_j^{M_{\nu}^{\epsilon}}\big(2^j(y-z_1),\dots,2^j(y-z_m)\big) \big| d\zzz\\
&\le 2^{j\nu n}\Big(\prod_{\mu=1}^{m}\ell(Q_{\mu})^{-n/p_{\mu}}\Big) \\
&\qq \q\times\int_{Q_1\times\cdots\times Q_{\nu}\times (\bbrn)^{m-\nu}} \big|K_j^{M_{\nu}^{\epsilon}}\big(2^j(y-z_1),\dots,2^j(y-z_{\nu}),z_{\nu+1},\dots,z_m\big) \big| d\zzz
\end{align*} with the obvious interpretations when $\nu=m$.
If $\nu=1$, then the integral in the preceding expression is bounded, using the Cauchy-Schwarz inequality,  by
\begin{align*}
&\ell(Q_1)^n\int_{Q_1}|Q_1|^{-1}\Big\Vert K_j^{M_{1}^{\epsilon}+\sum_{\mu\in J_m\setminus J_0}(n/p_{\mu}+\epsilon)}\big(2^j(y-z_1),z_2,\dots,z_m \big)\Big\Vert_{L^2(z_2,\dots,z_m)} dz_1\\
&\lesssim \ell(Q_1)^nh_{j,s}^{Q_1,0}(y)
\end{align*} where the inequality follows from the fact that 
$s>n/p-n/2+(m+1)\epsilon=M_1^{\epsilon}+\sum_{\mu\in J_m\setminus J_0}(n/p_{\mu}+\epsilon)$, which is due to \eqref{conditionons}.
When $\nu\ge 2$, we estimate the integral by
\begin{align*}
\Big(\prod_{\mu=1}^{\nu}\ell(Q_{\mu})^n \Big)\int_{Q_1}|Q_1|^{-1}\Big( \int_{(\bbrn)^{m-\nu}}  \Big\Vert K_j^{M_{\nu}^{\epsilon}}\big(2^j(y-z_1),z_2,\dots,z_m\big)\Big\Vert_{L^{\infty}(z_2,\dots,z_{\nu})}      dz_{\nu+1}\cdots dz_m\Big)dz_1
\end{align*}
and then, in view of Lemma \ref{imbedding}, we may replace the $L^{\infty}$ norm inside the integral by $L^2$ norm.
Now by applying the Cauchy-Schwarz inequality, the last expression is controlled by
\begin{align*}
&\Big(\prod_{\mu=1}^{\nu}\ell(Q_{\mu})^n \Big)\int_{Q_1}|Q_1|^{-1}  \Big\Vert K_j^{M_{\nu}^{\epsilon}+\sum_{\mu\in J_m\setminus J_0}(n/p_{\mu}+\epsilon)}\big(2^j(y-z_1),z_2,\dots,z_m\big)\Big\Vert_{L^{2}(z_2,\dots,z_{m})}      \; dz_1\\
&\lesssim \Big(\prod_{\mu=1}^{\nu}\ell(Q_{\mu})^n \Big)h_{j,s}^{Q_1,0}(y)
\end{align*}
as $s>M_{\nu}^{\epsilon}+\sum_{\mu\in J_m\setminus J_0}(n/p_{\mu}+\epsilon)$.
In both cases ( $\nu=1$ and $\nu\ge 2$ ), we have
\begin{align}\label{esthjq0}
\big| T_j\big(a_1,\dots,a_m\big)(y)\big|&\lesssim 2^{j\nu n} \Big(\prod_{\mu=1}^{\nu}\ell(Q_{\mu})^{-n/p_{\mu}+n}\Big)\Big( \prod_{\mu\in J_m\setminus J_0}\ell(Q_{\mu})^{-n/p_{\mu}}\Big)\\
&\qq\qq\qq\qq\times \frac{1}{\langle 2^j(y-c_{Q_1})\rangle^{N_{\nu}+\epsilon_{\nu}+n/2+\epsilon}} h_{j,\nu}^{Q_1,0}(y).\nonumber
\end{align}

Moreover, it follows from the vanishing moment condition of $a_{1}$ that
\begin{align}\label{tjay}
&\big| T_j\big(a_1,\dots,a_m\big)(y)\big| \nonumber\\
&\lesssim \big(2^j\ell(Q) \big)^L2^{jmn}\Big(\prod_{\mu=1}^{m}\ell(Q_{\mu})^{-n/p_{\mu}} \Big)\sum_{|\alpha|=L}\int_0^1\Big[ \int_{Q_1\times \cdots\times Q_m}\\
&\qq\qq\qq\qq\qq\times \big| \partial_{1}^{\alpha}K_j\big(2^jy_{c_{Q_1},z_1}^{\theta},2^j(y-z_2),\dots,2^j(y-z_m) \big)\big| d\zzz\Big]d\theta \nonumber
\end{align}
where we recall that $L$ is the number of vanishing moment of $a_1$, which is greater than $n/p-n$.
Similar to \eqref{ycqyz}, we have
$$|y-c_{Q_1}|\lesssim |y_{c_{Q_1},z_1}^{\theta}|$$
for $y\in (Q_1^*)^c$, $z_1\in Q_1$, and $0<\theta<1$. This deduces that 
\begin{align*}
&\langle 2^j(y-c_{Q_1})\rangle^{M_{\nu}^{\epsilon}}\int_{Q_1\times \cdots\times Q_m}\big| \partial_{1}^{\alpha}K_j\big(2^jy_{c_{Q_1},z_1}^{\theta},2^j(y-z_2),\dots,2^j(y-z_m) \big)\big| d\zzz\\
&\lesssim 2^{-j(m-\nu)n}\int_{Q_1\times \cdots\times Q_{\nu}\times (\bbrn)^{m-\nu}}    \big| \big(\partial_{1}^{\alpha}K_j\big)^{M_{\nu}^{\epsilon}}\big(2^jy_{c_{Q_1},z_1}^{\theta},2^j(y-z_2),\dots,2^j(y-z_{\nu}),z_{\nu+1},\dots,z_m \big)\big| d\zzz
\end{align*}
with the obvious interpretation when $\nu=m$.
Now using the arguments that led to \eqref{esthjq0}, we can bound the preceding expression by
\begin{align*}
&2^{-j(m-\nu)n}\Big( \prod_{\mu=1}^{\nu}\ell(Q_{\mu})^n \Big)\int_{Q_1}|Q_1|^{-1} \Big\Vert  K_j^{M_{\nu}^{\epsilon}+\sum_{\mu\in J_m\setminus J_0}(n/p_{\mu}+\epsilon)}\big(2^jy_{c_{Q_1},z_1}^{\theta},z_2,\dots,z_m\big)    \Big\Vert_{L^2(z_2,\dots,z_m)}dz_1\\
&\lesssim 2^{-j(m-\nu)n}\Big( \prod_{\mu=1}^{\nu}\ell(Q_{\mu})^n \Big) \int_{Q_1}|Q_1|^{-1} \Big\Vert  K_j^{s}\big(2^jy_{c_{Q_1},z_1}^{\theta},z_2,\dots,z_m\big)    \Big\Vert_{L^2(z_2,\dots,z_m)}dz_1
\end{align*}
as $s>M_{\nu}^{\epsilon}+\sum_{\mu\in J_m\setminus J_0}(n/p_{\mu}+\epsilon)$.
Combining all together, we finally obtain
\begin{align*}
\big| T_j\big(a_1,\dots,a_m\big)(y)\big|&\lesssim 2^{j\nu n} \Big(\prod_{\mu=1}^{\nu}\ell(Q_{\mu})^{-n/p_{\mu}+n}\Big)\Big( \prod_{\mu\in J_m\setminus J_0}\ell(Q_{\mu})^{-n/p_{\mu}}\Big)\\
&\qq\qq\qq\qq\times \frac{1}{\langle 2^j(y-c_{Q_1})\rangle^{N_{\nu}+\epsilon_{\nu}+n/2+\epsilon}} h_{j,s}^{Q_1,L}(y).\nonumber
\end{align*}
This, together with the estimate \eqref{esthjq0}, completes the proof of \eqref{keyestimate1}.

\subsection{The second assertion}
Suppose that
$$x\in  \Big( \bigcap_{\mu \in J_0\setminus I_0}(Q_{\mu}^{**})^c\cap B_{\mu,l}\Big)\medcap \Big(\bigcap_{\mu\in I_0}(Q_{\mu}^{**})^c\cap (B_{\mu,l})^c  \Big)     \medcap \Big( \bigcap_{\mu\in J_m\setminus J_0}Q_{\mu}^{**}\;\Big).$$
Since $I_0\not=\emptyset$, we can choose $\nu_0\in I_0$ where we set $\nu_0=1$ if $I_0$ contains $1$.
We first note that for any $\mu\in J_0\setminus I_0$, 
$$|x-c_{Q_{\mu}}|\le |x-c_{Q_{\nu_0}}|$$
since $x\in B_{\mu,l}\cap B_{\nu_0,l}^c$.
Moreover, if $\mu\in I_0$, $x\in (Q_{\mu}^{**})^c\cap (B_{\mu,l})^c$, and $|x-y|\le 2^{-l}$, then
we have 
$$|x-c_{Q_{\mu}}|\lesssim |y-c_{Q_{\mu}}|.$$
Therefore, it follows that
\begin{align}\label{keykeyest1}
&\big\langle 2^j(x-c_{Q_{1}})\big\rangle^{n/p_1-n/2+\epsilon}\Big( \prod_{\mu=2}^{\nu}\big\langle 2^j(x-c_{Q_{\mu}}) \big\rangle^{n/p_{\mu}+\epsilon}\Big) \big| \phi_l\ast T_j\big(a_1,\dots,a_m\big)(x)\big|\nonumber\\
&\lesssim 2^{ln}\int_{|x-y|\le 2^{-l}}   \big\langle 2^j(y-c_{Q_{\nu_0}})\big\rangle^{-n/2+\sum_{\mu\in J_0\setminus I_0}(n/p_{\mu}+\epsilon)}\\
&\qq\qq\qq \times \Big(\prod_{\mu\in I_0}\big\langle 2^j(y-c_{Q_{\mu}}) \big\rangle^{n/p_{\mu}+\epsilon} \Big) \big| T_j\big(a_1,\dots,a_m\big)(y)\big|\chi_{\cap_{\mu\in I_0}(Q_{\mu}^{*})^c}(y) dy \nonumber
\end{align}
as the conditions $x\in (Q_{\mu}^{**})^c\cap (B_{\mu,l})^c$ and  $|x-y|\le 2^{l}$ imply $y\in (Q_{\mu}^{*})^c$. 
In order to proceed further we claim that the following lemma holds.
\begin{lemma}\label{lem keykeyest2}
For $y\in \bigcap_{\mu\in I_0}(Q_{\mu}^*)^c$, 
\begin{align}\label{keykeyest2}
& \big\langle 2^j(y-c_{Q_{\nu_0}})\big\rangle^{ -n/2+\sum_{\mu\in J_0\setminus I_0}(n/p_{\mu}+\epsilon)}\Big(\prod_{\mu\in I_0}\big\langle 2^j(y-c_{Q_{\mu}}) \big\rangle^{n/p_{\mu}+\epsilon} \Big)\big| T_j\big(a_1,\dots,a_m\big)(y)\big| \nonumber\\
 &\lesssim 2^{j\nu n} \Big(\prod_{\mu=1}^{m}\ell(Q_{\mu})^{-n/p_{\mu}} \Big)\Big(\prod_{\mu=1}^{\nu}\ell(Q_{\mu})^n \Big)\min\big\{ h_{j,s}^{Q_1,0}(y),h_{j,s}^{Q_1,L}(y) \big\}
\end{align}
where $h_{j,s}^{Q_1,0}$ and $h_{j,s}^{Q_1,L}$ are defined as in \eqref{hjs1def} and \eqref{hjs1ldef}.
\end{lemma}
The proof of \eqref{keykeyest2} will be given in the last part of this section.
For now we take \eqref{keykeyest2} for granted.
Then the left-hand side of \eqref{keykeyest1} is bounded by a constant multiple of
$$2^{j\nu n}\Big(\prod_{\mu=1}^{m}\ell(Q_{\mu})^{-n/p_{\mu}} \Big)\Big(\prod_{\mu=1}^{\nu}\ell(Q_{\mu})^n \Big)\min\big\{ \mathcal{M}h_{j,s}^{Q_1,0}(x),\mathcal{M}h_{j,s}^{Q_1,L}(x) \big\}, $$
which, together with \eqref{lpdecomp}, implies that
\begin{align*}
&\big| \phi_l\ast T_{\sigma}\big(a_1,\dots,a_m\big)(x)\big|\\
&\lesssim \sum_{j\in\bbz}2^{j\nu n}\Big(\prod_{\mu\in J_m\setminus J_0} \ell(Q_{\mu})^{-n/p_{\mu}}\chi_{Q_{\mu}^{**}}(x)\Big)\Big(\prod_{\mu=2}^{\nu}\frac{\ell(Q_{\mu})^{-(n/p_{\mu}-n)}}{\langle {2^j(x-c_{Q_{\mu}})}\rangle^{n/p_{\mu}+\epsilon}} \chi_{(Q_{\mu}^{**})^c}(x)\Big)\\
&\qq\qq\times \frac{\ell(Q_1)^{-(n/p_1-n)}}{\langle 2^j(x-c_{Q_1})\rangle^{n/p_1-n/2+\epsilon}}\chi_{(Q_{1}^{**})^c}(x)\min\big\{ \mathcal{M}h_{j,s}^{Q_1,0}(x),\mathcal{M}h_{j,s}^{Q_1,L}(x) \big\}\\
&\le \Big(\prod_{\mu\in J_m\setminus J_0} \ell(Q_{\mu})^{-n/p_{\mu}}\chi_{Q_{\mu}^{**}}(x)\Big) \Big(\prod_{\mu=2}^{\nu}\frac{1}{| {x-c_{Q_{\mu}}}|^{n/p_{\mu}+\epsilon}} \chi_{(Q_{\mu}^{**})^c}(x)\Big) \ell(Q_1)^{-\sum_{\mu=1}^{\nu}(n/p_{\mu}-n)}\\
&\qq\qq\times \frac{\chi_{(Q_{1}^{**})^c}(x)}{|x-c_{Q_1}|^{n/p_1-n/2+\epsilon}} \sum_{j\in\bbz}2^{j\nu n}2^{jn/2}2^{-j(\sum_{\mu=1}^{\nu}(n/p_{\mu}+\epsilon))}\min\big\{ \mathcal{M}h_{j,s}^{Q_1,0}(x),\mathcal{M}h_{j,s}^{Q_1,L}(x) \big\}\\
&\le \Big(\prod_{\mu\in J_m\setminus J_0} \ell(Q_{\mu})^{-n/p_{\mu}}\chi_{Q_{\mu}^{**}}(x)\Big) \Big(\prod_{\mu=2}^{\nu}\frac{\ell(Q_{\mu})^{\epsilon}}{| {x-c_{Q_{\mu}})}|^{n/p_{\mu}+\epsilon}} \chi_{(Q_{\mu}^{**})^c}(x)\Big) \ell(Q_1)^{-\sum_{\mu=1}^{\nu}(n/p_{\mu}-n+\epsilon)}\\
&\qq\times \frac{\ell(Q_1)^{\epsilon}}{|x-c_{Q_1}|^{n/p_1-n/2+\epsilon}} \chi_{(Q_{1}^{**})^c}(x)\sum_{j\in\bbz}2^{jn/2}2^{-j(\sum_{\mu=1}^{\nu}(n/p_{\mu}-n+\epsilon))}\min\big\{ \mathcal{M}h_{j,s}^{Q_1,0}(x),\mathcal{M}h_{j,s}^{Q_1,L}(x) \big\}.
\end{align*}
Here, we used the assumption that $\ell(Q_1)\le \ell(Q_{\mu})$, $\mu=2,\dots,\nu$, in the last two inequalities.
Finally, the inequality \eqref{keylemma2assertion} follows from taking
\begin{align*}
\La_1^{J_0}(x)&:=\big( \LL_s^2[\sigma]\big)^{-1}\ell(Q_1)^{-\sum_{\mu=1}^{\nu}(n/p_{\mu}-n+\epsilon)}\frac{\ell(Q_1)^{\epsilon}}{|x-c_{Q_1}|^{n/p_1-n/2+\epsilon}}\chi_{(Q_{1}^{**})^c}(x) \\
&\qq\times \sum_{j\in\bbz}2^{jn/2}2^{-j(\sum_{\mu=1}^{\nu}(n/p_{\mu}-n+\epsilon))}\min\big\{ \mathcal{M}h_{j,s}^{Q_1,0}(x),\mathcal{M}h_{j,s}^{Q_1,L}(x) \big\}\\
\La_{\mu}^{J_0}(x)&:=\begin{cases}
\frac{\ell(Q_{\mu})^{\epsilon}}{|x-c_{Q_{\mu}}|^{n/p_{\mu}+\epsilon}}\chi_{(Q_{\mu}^{**})^c}(x), \qq& \mu\in J_0\setminus \{1\}\\
\ell(Q_{\mu})^{-n/p_{\mu}}\chi_{Q_{\mu}^{**}}(x), \qq& \mu\in J_m\setminus J_0
\end{cases}.
\end{align*}
Moreover, using H\"older's inequality with $(1/p_1-1/2)+1/2=1/p_1$,
\begin{align*}
\big\Vert \La_1^{J_0}\big\Vert_{L^{p_1}(\bbrn)}&\le \big( \LL_s^2[\sigma]\big)^{-1}\ell(Q_1)^{-\sum_{\mu=1}^{\nu}(n/p_{\mu}-n+\epsilon)}\bigg\Vert \frac{\ell(Q_1)^{\epsilon}}{|\cdot-c_{Q_1}|^{n/p_1-n/2+\epsilon}}\bigg\Vert_{L^{\frac{1}{1/p_1-1/2}}((Q_1^{**})^c)}\\
&\qq\times \Big\Vert \sum_{j\in\bbz}2^{jn/2}2^{-j(\sum_{\mu=1}^{\nu}(n/p_{\mu}-n+\epsilon))}\min\big\{ \mathcal{M}h_{j,s}^{Q_1,0},\mathcal{M}h_{j,s}^{Q_1,L} \big\} \Big\Vert_{L^2(\bbrn)}.
\end{align*}
The $L^{\frac{1}{1/p_1-1/2}}$ norm is controlled by a constant and the $L^2$ norm is less than
\begin{align*}
&\sum_{j\in\bbz}2^{jn/2}2^{-j(\sum_{\mu=1}^{\nu}(n/p_{\mu}-n+\epsilon))}\min\Big\{ \big\Vert \mathcal{M}h_{j,s}^{Q_1,0}\big\Vert_{L^2(\bbrn)}, \big\Vert \mathcal{M}h_{j,s}^{Q_1,L}\big\Vert_{L^2(\bbrn)} \Big\}\\
&\lesssim \LL_s^2[\sigma]\sum_{j\in\bbz}2^{-j(\sum_{\mu=1}^{\nu}(n/p_{\mu}-n+\epsilon))}\min\Big\{1,\big(2^j\ell(Q_1)\big)^L \Big\}\\
&\lesssim \LL_s^2[\sigma]\ell(Q_1)^{\sum_{\mu=1}^{\nu}(n/p_{\mu}-n+\epsilon)}
\end{align*}
for $L>n/p-n\ge n/p-mn+m\epsilon$,
where we applied the $L^2$ boundedness of $\mathcal{M}$ and the estimates \eqref{hjl2norms} in the first inequality.
This completes the proof of 
$$\big\Vert \La_1^{J_0}\big\Vert_{L^{p_1}(\bbrn)}\lesssim 1.$$
It is also clear that $$\big\Vert \La_{\mu}^{J_0}\big\Vert_{L^{p_1}(\bbrn)}\lesssim 1, \qq\mu=2,\dots,m.$$

\hfill

\subsection{Proof of Lemma \ref{lem keykeyest2} }

Now let us prove the claim \eqref{keykeyest2}.
Suppose that $y\in \bigcap_{\mu \in I_0}(Q_{\mu}^*)^c$. For simplicity, we use the notations
$$s_{\mu}:=\begin{cases}
n/p_{\nu_0}-n/2+\epsilon+\sum_{\eta \in J_0\setminus I_0}(n/p_{\eta}+\epsilon), \qq& \mu=\nu_0\\
n/p_{\mu}+\epsilon,\qq & \mu\not= \nu_0,~ \mu\in I_0
\end{cases}$$
and
$$S_{I_0}:=\sum_{\mu\in I_0}s_{\mu} = \Big( \sum_{\mu \in J_0} \frac{n}{p_\mu} \Big) - \frac{n}{2} + \nu\epsilon.$$
Then the left-hand side of \eqref{keykeyest2} can be written as
\begin{equation}\label{mainestterm}
\Big(\prod_{\mu\in I_0}{\big\langle 2^j(y-c_{Q_{\mu}})\big\rangle^{s_{\mu}}}\Big)\big| T_j\big(a_1,\dots,a_m\big)(y)\big|.
\end{equation}

\hfill

We first observe that for $y\in (Q_{\mu}^*)^c$ and $z_{\mu}\in Q_{\mu}$,
\begin{equation}\label{ycqmuyzmu}
|y-c_{Q_{\mu}}|\lesssim |y-z_{\mu}|
\end{equation}
and this implies that \eqref{mainestterm} is controlled by a constant times
\begin{align*}
2^{jmn}\Big(\prod_{\mu=1}^{m}\ell(Q_{\mu})^{-n/p_{\mu}} \Big)\int_{Q_1\times\cdots\times Q_m}\big| K_j^{S_{I_0}}\big(2^j(y-z_1),\dots,2^j(y-z_m)\big)\big| d\zzz.
\end{align*}
where $K_j^r$ is defined as in \eqref{kjrdef}.
Now we apply a change of variables to bound the above expression by
\begin{align}\label{arginfirst}
& 2^{j\nu n}\Big(\prod_{\mu=1}^{m}\ell(Q_{\mu})^{-n/p_{\mu}} \Big)\int_{Q_1\times\cdots\times Q_{\nu}\times (\bbrn)^{m-\nu}}    \big| K_j^{S_{I_0}}\big(2^j(y-z_1),\dots,2^j(y-z_\nu),z_{\nu+1},\dots,z_m\big)\big| d\zzz \nonumber\\
&\le 2^{j\nu n}\Big(\prod_{\mu=1}^{m}\ell(Q_{\mu})^{-n/p_{\mu}} \Big)\Big(\prod_{\mu=1}^{\nu}\ell(Q_{\mu})^n \Big)\int_{Q_1}|Q_1|^{-1}\\
&\qq\qq\qq\times \Big( \int_{(\bbrn)^{m-\nu}}  \big\Vert K_j^{S_{I_0}}\big(2^j(y-z_1),z_2,\dots,z_m\big)\big\Vert_{L^{\infty}(z_2,\dots,z_{\nu})}   dz_{\nu+1}\cdots dz_m \Big)\;dz_1 \nonumber
\end{align}
with the obvious interpretations when $\nu=1$ or $\nu= m$ as in the proof of \eqref{esthjq0}.
In view of Lemma \ref{imbedding}, we may replace the $L^{\infty}$ norm inside the integral by $L^2$ norm and then apply the Cauchy-Schwarz inequality to dominate the above expression by
\begin{align*}
&2^{j\nu n}\Big(\prod_{\mu=1}^{m}\ell(Q_{\mu})^{-n/p_{\mu}} \Big)\Big(\prod_{\mu=1}^{\nu}\ell(Q_{\mu})^n \Big)\int_{Q_1}|Q_1|^{-1}\\
&\qq\qq\times \Big\Vert K_j^{S_{I_0}+\sum_{\mu\in J_m\setminus J_0}(n/p_{\mu}+\epsilon)}\big(2^j(y-z_1),z_2,\dots,z_m\big)\Big\Vert_{L^2(z_2,\dots,z_m)} dz_1\\
&\le 2^{j\nu n}\Big(\prod_{\mu=1}^{m}\ell(Q_{\mu})^{-n/p_{\mu}} \Big)\Big(\prod_{\mu=1}^{\nu}\ell(Q_{\nu})^n \Big)h_{j,s}^{Q_1,0}(y)
\end{align*}
as \begin{equation}\label{conditionss}
s>n/p-n/2+m\epsilon=S_{I_0}+\sum_{\mu\in J_m\setminus J_0}(n/p_{\mu}+\epsilon).
\end{equation}

Moreover, using \eqref{tjay} and \eqref{ycqmuyzmu}, the term \eqref{mainestterm} is no more than a constant multiple of
 \begin{align}\label{smutjest}
 &\big(2^j\ell(Q_1) \big)^L 2^{jmn}\Big(\prod_{\mu=1}^{m}\ell(Q_{\mu})^{-n/p_{\mu}} \Big)  \sum_{|\alpha|=L}\int_0^1\Big[   \int_{Q_1\times\cdots\times Q_m}\\
 &\qq\qq\qq\qq\times \big| \big(\partial_1^{\alpha}K_j\big)^{S_{I_0}}\big(2^jy_{c_{Q_1},z_1}^{\theta},2^j(y-z_2),\dots,2^j(y-z_m)\big)\big| d\zzz \Big] d\theta \nonumber
 \end{align}
 where  $\big(\partial_1^{\alpha}K_j\big)^r$ is defined as in \eqref{kjrdef}.
In the above estimate, if $I_0$ contains $1$, then we also applied the fact that  for $y\in (Q_1^*)^c$, $z_1\in Q_1$, and $0<\theta<1$,
 $$|y-c_{Q_1}|\lesssim \big| y_{c_{Q_1},z_1}^{\theta}\big|.$$
Then, similar to \eqref{arginfirst}, we can bound  \eqref{smutjest} by
\begin{align*}
&2^{j\nu n}\Big(\prod_{\mu=1}^{m}\ell(Q_{\mu})^{-n/p_{\mu}} \Big)\Big(\prod_{\mu=1}^{\nu}\ell(Q_{\mu})^n \Big)     \big( 2^j\ell(Q_1)\big)^L\sum_{|\alpha|=L}\int_0^1 \Big[ \int_{Q_1}|Q_1|^{-1}\\
&\qq\qq\times \Big(\int_{(\bbrn)^{m-\nu}} \big\Vert \big( \partial_1^{\alpha}K_j\big)^{S_{I_0}}\big(2^j y_{c_{Q_1},z_1}^{\theta},z_2,\dots,z_m\big)    \big\Vert_{L^{\infty}(z_2,\dots,z_{\nu})}dz_{\nu+1}\cdots dz_m \Big)        dz_1\Big]\; d\theta
\end{align*}
with the usual modification when $\nu=1$ or $\nu=m$.
Clearly, the $L^{\infty}$ norm inside the integral can be replaced, thanks to Lemma \ref{imbedding}, by 
$$\big\Vert K_j^{S_{I_0}}\big(2^jy_{c_{Q_1},z_1}^{\theta},z_2,\dots,z_m \big)\big\Vert_{L^2(z_2,\dots,z_{\nu})}$$
and thus we finally arrive at the inequality, via the Cauchy-Schwarz inequality,
\begin{align*}
\eqref{mainestterm}&\lesssim 2^{j\nu n}\Big(\prod_{\mu=1}^{m}\ell(Q_{\mu})^{-n/p_{\mu}} \Big)\Big(\prod_{\mu=1}^{\nu}\ell(Q_{\mu})^n \Big) h_{j,s}^{Q_1,L}(y).
\end{align*}
with the aid of \eqref{conditionss}.
This ends the proof of Lemma \ref{lem keykeyest2}.

\section*{Acknowledgement}

{We would like to thank the anonymous referees for the careful reading and their helpful comments.}

\end{document}